\documentclass[11pt]{article}

\usepackage{amsmath}
\usepackage{amssymb}
\usepackage{amsthm}
\usepackage{latexsym}
\usepackage{color}
\usepackage{graphicx}
\usepackage{color} 
\usepackage[T1]{fontenc}
\usepackage[english]{babel}

\DeclareSymbolFont{calletters}{OMS}{cmsy}{m}{n}
\DeclareSymbolFontAlphabet{\mathcal}{calletters}
\def\be{\begin{eqnarray}}
\def\ee{\end{eqnarray}}

\def\b*{\begin{eqnarray*}}
\def\e*{\end{eqnarray*}}

\newtheorem{Theorem}{Theorem}[part]
\newtheorem{Definition}{Definition}[part]
\newtheorem{Proposition}{Proposition}[part]

\newtheorem{Assumption}{Assumption}[part]
\newtheorem{Lemma}{Lemma}[part]

\newtheorem{Remark}{Remark}[part]

\makeatletter \@addtoreset{equation}{section}

\@addtoreset{Definition}{section}

\@addtoreset{Theorem}{section}

\@addtoreset{Proposition}{section}

\@addtoreset{Property}{section}

\@addtoreset{Assumption}{section}

\@addtoreset{Corollary}{section}

\@addtoreset{Lemma}{section}

\@addtoreset{Remark}{section}

\@addtoreset{Example}{section}

\newcommand{\No}[1]{\left\|#1\right\|}     % Norm
\newcommand{\abs}[1]{\left|#1\right|}     % absolute value

\def\Fh{{\hat F}}

\def \Inf{\displaystyle\inf}
\def \Sup{\displaystyle\sup}

\def\einf{{\rm ess \, inf}}
\def\esup{{\rm ess \, sup}}
\def\trace{{\rm Tr}}

\def\={\;=\;}
\def\.{\;.}

\def\eps{\varepsilon}

\def\reff#1{{\rm(\ref{#1})}}

\def\1{{\bf 1}}
\def \ep{\hbox{ }\hfill{ ${\cal t}$~\hspace{-5.1mm}~${\cal u}$   } }

\def \proof{{\noindent \bf Proof. }}
\def \ep{\hbox{ }\hfill$\Box$}

 \def\normeL2#1{\left\|{#1}\right\|_{L^2}}
 \title{Second Order Backward Stochastic Differential Equations under Monotonicity Condition \footnote{Research partly supported by the Chair {\it Financial Risks} of the {\it Risk Foundation} sponsored by Soci\'et\'e G\'en\'erale, the Chair {\it Derivatives of the Future} sponsored by the {F\'ed\'eration Bancaire Fran\c{c}aise}, and the Chair {\it Finance and Sustainable Development} sponsored by EDF and Calyon.}}
\author{ Dylan {\sc Possama\"{i}}\footnote{CMAP, Ecole Polytechnique, Paris, dylan.possamai@polytechnique.edu.} }          
      \date{\today}

\begin{document}

\maketitle

\begin{abstract}
In a recent paper, Soner, Touzi and Zhang \cite{stz} have introduced a notion of second order backward stochastic differential equations (2BSDEs for short), which are naturally linked to a class of fully non-linear PDEs. They proved existence and uniqueness for a generator which is uniformly Lipschitz in the variables $y$ and $z$. The aim of this paper is to extend these results to the case of a generator satisfying a monotonicity condition in $y$. More precisely, we prove existence and uniqueness for 2BSDEs with a generator which is Lipschitz in $z$ and uniformly continuous with linear growth in $y$. Moreover, we emphasize throughout the paper the major difficulties and differences due to the 2BSDE framework.
\vspace{10mm}

\vspace{5mm}

\end{abstract}
\newpage

\section{Introduction}

Backward stochastic differential equations (BSDEs for short) appeared in Bismut \cite{bis} in the linear case, and then have been widely studied since the seminal paper of Pardoux and Peng \cite{pardpeng}. Their range of applications includes notably probabilistic numerical methods for partial differential equations, stochastic control, stochastic differential games, theoretical economics and financial mathematics.

\vspace{0.5em}
On a filtered probability space $(\Omega,\mathcal F,\left\{\mathcal F_t\right\}_{0\leq t\leq T},\mathbb P)$ generated by an $\mathbb R^d$-valued Brownian motion $B$, a solution to a BSDE consists on finding a pair of progressively measurable processes $(Y,Z)$ such that
$$Y_t=\xi +\int_t^Tf_s(Y_s,Z_s)ds-\int_t^TZ_sdB_s,\text{ } t\in [0,T], \text{ }\mathbb P-a.s.$$
where $f$ (also called the driver) is a progressively measurable function and $\xi$ is an $\mathcal F_T$-measurable random variable.

\vspace{0.5em}
Pardoux and Peng \cite{pardpeng} proved existence and uniqueness of the above BSDE provided that the function $f$ is uniformly Lipschitz in $y$ and $z$ and that $\xi$ and $f_s(0,0)$ are square integrable. Then, in \cite{pardpeng2}, they proved that if the randomness in $f$ and $\xi$ is induced by the current value of a state process defined by a forward stochastic differential equation, then the solution to the BSDE could be linked to the solution of a semilinear PDE by means of a generalized Feynman-Kac formula. Since their pioneering work, many efforts have been made to relax the assumptions on the driver $f$. For instance, Lepeltier and San Martin \cite{lsm} have proved the existence of a solution when $f$ is only continuous in $(y,z)$ with linear growth.

\vspace{0.8em}
More recently, motivated by applications in financial mathematics and probabilistic numerical methods for PDEs (see \cite{ftw}), Cheredito, Soner, Touzi and Victoir \cite{cstv} introduced the notion of second order BSDEs (2BSDEs), which are connected to the larger class of fully nonlinear PDEs. Then, Soner, Touzi and Zhang \cite{stz} provided a complete theory of existence and uniqueness for 2BSDEs under uniform Lipschitz conditions similar to those of Pardoux and Peng. Their key idea was to reinforce the condition that the 2BSDE must hold $\mathbb P-a.s.$ for every probability measure $\mathbb P$ in a non-dominated class of mutually singular measures (see Section \ref{section.1} for precise definitions). Let us describe the intuition behind their formulation.

\vspace{0.8em}
Suppose that we want to study the following fully non-linear PDE
\begin{equation}
- \frac{\partial u}{\partial t}-h\left(t,x,u(t,x),Du(t,x),D^2u(t,x)\right)=0,\qquad u(T,x)=g(x).
\label{eq:1}
\end{equation}

\vspace{0.8em}
If the function $\gamma\mapsto h(t,x,r,p,\gamma)$ is assumed to be convex, then it is equal to its double Fenchel-Legendre transform, and if we denote its Fenchel-Legendre transform by $f$, we have
\begin{equation}
h(t,r,p,\gamma)=\underset{a\geq 0}{\sup}\left\{\frac12a\gamma-f(t,x,r,p,a)\right\}
\label{eq:2}
\end{equation}

\vspace{0.8em}
Then, from \reff{eq:2}, we expect, at least formally, that the solution $u$ of \reff{eq:1} is going to verify
$$u(t,x)=\underset{a\geq 0}{\sup}\ u^a(t,x),$$
where $u^a$ is defined as the solution of the following semi-linear PDE
\begin{equation}
- \frac{\partial u^a}{\partial t}-\frac12aD^2u^a(t,x)+f\left(t,x,u^a(t,x),Du^a(t,x),a\right)=0,\qquad u^a(T,x)=g(x).
\label{eq:3}
\end{equation}

\vspace{0.8em}
Since $u^a$ is linked to a classical BSDE, the 2BSDE associated to $u$ should correspond (in some sense) to the supremum of the family of BSDEs indexed by $a$. Furthermore, changing the process $a$ can be achieved by changing the probability measure under which the BSDE is written. In these respects, the 2BSDE theory shares deep links with the theory of quasi-sure stochastic analysis of Denis and Martini \cite{denis} and the theory of G-expectation of Peng \cite{peng}.

\vspace{0.8em}
In addition to providing a probability representation for solutions of fully non-linear PDEs, the 2BSDE theory also has many applications in mathematical finance, especially within the context of markets with volatility uncertainty. Hence, Soner, Touzi and Zhang \cite{stz3} solved the superhedging problem under volatility uncertainty within this quasi-sure framework (following the original appoach of \cite{denis}). More recently, Matoussi, Possama\"{i} and Zhou \cite{mpz} proved that the solution of the utility maximization problem for an investor in an incomplete market with volatility uncertainty could be expressed in terms of a second order BSDE. The same authors also used this theory in \cite{mpz2} to give a superhedging price for American options under volatility uncertainty.

\vspace{0.8em}
Our aim in this paper is to relax the Lipschitz-type hypotheses of \cite{stz} on the driver of the 2BSDE to prove an existence and uniqueness result. In Section \ref{section.1}, inspired by Pardoux \cite{pard}, we study 2BSDEs with a driver which is Lipschitz in some sense in $z$, uniformly continuous with linear growth in $y$ and satisfies a monotonicity condition. We then prove existence and uniqueness and highlight one of the main difficulties when dealing with 2BSDEs. Indeed, the main tool in the proof of existence is to use monotonic approximations (as in \cite{lsm}). However, since we are working under a family of non-dominated probability measures, the monotone or dominated convergence theorem may fail, which in turn raises subtle technical difficulties in the proofs.

\section{Preliminaries} \label{section.1}

Let $\Omega:=\left\{\omega\in C([0,T],\mathbb R^d):\omega_0=0\right\}$ be the canonical space equipped with the uniform norm $\No{\omega}_{\infty}:=\sup_{0\leq t\leq T}|\omega_t|$, $B$ the canonical process, $\mathbb P_0$ the Wiener measure, $\mathbb F:=\left\{\mathcal F_t\right\}_{0\leq t\leq T}$ the filtration generated by $B$, and $\mathbb F^+:=\left\{\mathcal F_t^+\right\}_{0\leq t\leq T}$ the right limit of $\mathbb F$. We first recall the notations introduced in \cite{stz}.

\subsection{The Local Martingale Measures}

We will say that a probability measure $\mathbb P$ is a local martingale measure if the canonical process $B$ is a local martingale under $\mathbb P$. By Karandikar \cite{kar}, we know that we can give pathwise definitions of the quadratic variation $\left<B\right>_t$ and its density $\widehat a_t$.

\vspace{0.8em}
Let $\overline{\mathcal P}_W$ denote the set of all local martingale measures $\mathbb P$ such that
\begin{equation}
\left<B\right>_t \text{ is absolutely continuous in $t$ and $\widehat a$ takes values in $\mathbb S_d^{>0}$, }\mathbb P-a.s.
\label{eq:}
\end{equation}
where $\mathbb S_d^{>0}$ denotes the space of all $d\times d$ real valued positive definite matrices.

\vspace{0.8em}
We recall from \cite{stz}, the class $\overline{\mathcal P}_S\subset\overline{\mathcal P}_W$ consisting of all probability measures
\begin{equation}\label{proba.p}
\mathbb P^\alpha:=\mathbb P_0\circ (X^\alpha)^{-1} \text{ where }X_t^\alpha:=\int_0^t\alpha_s^{1/2}dB_s,\text{ }t\in [0,1],\text{ }\mathbb P_0-a.s.,
\end{equation}

and we concentrate on the subclass $\widetilde{\mathcal P}_S\subset\overline{\mathcal P}_S$ defined by
\begin{equation}\label{ps}
\widetilde{\mathcal P}_S:=\left\{\mathbb P^\alpha\in\overline{\mathcal P}_S,\ \underline a\leq\alpha\leq\bar a, \ \mathbb P_0-a.s.\right\},
\end{equation}
for fixed matrices $\underline a$ and $\bar a$ in $\mathbb S_d^{>0}$. We recall from \cite{stz2} that every $\mathbb P\in \overline{\mathcal P}_S$ (and thus in $\widetilde{\mathcal P}_S$) satisfies the Blumenthal zero-one law and the martingale representation property.

\vspace{0.8em}
We finish with a definition. Let $\mathcal P_0\subset\overline{\mathcal P}_W$.

\vspace{0.8em}
\begin{Definition}
We say that a property holds $\mathcal P_0$-quasi surely ($\mathcal P_0-q.s.$ for short) if it holds $\mathbb P-a.s.$ for all $\mathbb P\in\mathcal P_0$.
\end{Definition}

\subsection{The non-linear Generator}
We consider a map $H_t(\omega,y,z,\gamma):[0,T]\times\Omega\times\mathbb{R}\times\mathbb{R}^d\times D_H\rightarrow \mathbb{R}$, where $D_H \subset \mathbb{R}^{d\times d}$ is a given subset containing $0$. 

\vspace{0.8em}
Following the PDE intuition outlined in the Introduction, we define the corresponding conjugate of $H$ w.r.t. $\gamma$ by
\begin{align*}
&F_t(\omega,y,z,a):=\underset{\gamma \in D_H}{\Sup}\left\{\frac12\trace(a\gamma)-H_t(\omega,y,z,\gamma)\right\} \text{ for } a \in \mathbb S_d^{>0}\\[0.3em]
&\widehat{F}_t(y,z):=F_t(y,z,\widehat{a}_t) \text{ and } \widehat{F}_t^0:=\widehat{F}_t(0,0).
\end{align*}

We denote by $D_{F_t(y,z)}:=\left\{a, \text{ }F_t(\omega,y,z,a)<\infty\right\}$ the domain of $F$ in $a$ for a fixed $(t,\omega,y,z)$.

\vspace{0.8em}
As in \cite{stz} we fix a constant $\kappa \in (1,2]$ and restrict the probability measures in

\begin{Definition}\label{def}
$\mathcal{P}_H^\kappa$ consists of all $\mathbb P \in \widetilde{\mathcal{P}}_S$ such that
$$ \mathbb{E}^{\mathbb{P}}\left[\left(\int_0^T\abs{\widehat{F}_t^0}^\kappa dt\right)^{\frac2\kappa}\right]<+\infty.$$
\end{Definition}

\vspace{0.8em}
It is clear that $\mathcal{P}^\kappa_H$ is decreasing in $\kappa$, and $\widehat{a}_t \in D_{F_t}$, $dt\times d\mathbb{P}-a.s.$ for all $\mathbb{P} \in \mathcal {P}_H^\kappa$. We will also denote $\overline{\mathcal P}^\kappa_H$ the closure for the weak topology of $\mathcal P^\kappa_H$.

\vspace{0.8em}
\begin{Remark}
Unlike in \cite{stz}, we assume that the bounds on the density of the quadratic variation $\widehat a$ are uniform with respect to the underlying probability measure. In particular, this ensures that the family $\mathcal P_H^\kappa$ is weakly relatively compact and that $\overline{\mathcal P}^\kappa_H$ is weakly compact.
\end{Remark}

\vspace{0.8em}
We now state our main assumptions on the function $F$ which will be our main interest in the sequel
\begin{Assumption} \label{assump.h}
\begin{itemize}
\item[\rm{(i)}] The domain $D_{F_t(y,z)}=D_{F_t}$ is independent of $(\omega,y,z)$.
\item[\rm{(ii)}] For fixed $(y,z,a)$, $F$ is $\mathbb{F}$-progressively measurable.
\item[\rm{(iii)}] We have the following uniform Lipschitz-type property
$$\forall (y,z,z^{'},a,t,\omega), \text{ } \abs{F_t(\omega,y,z,a)- F_t(\omega,y,z^{'},a)}\leq C\abs{a^{1/2}(z-z^{'})}.$$
\item[\rm{(iv)}] $F$ is uniformly continuous in $\omega$ for the $||\cdot||_\infty$ norm.
\item[\rm{(v)}] $F$ is uniformly continuous in $y$, uniformly in $(z,t,\omega,a)$, and has the following growth property
$$\exists C>0 \text{ s.t. } \forall (t,y,a,\omega), \text{ }\abs{F_t(\omega,y,0,a)}\leq \abs{F_t(\omega,0,0,a)}+C(1+\abs{y}).$$
\item[\rm{(vi)}] We have the following monotonicity condition. There exists $\mu>0$ such that
$$ (y_1-y_2)(F_t(\omega,y_1,z,a)-F_t(\omega,y_2,z,a))\leq \mu\abs{y_1-y_2}^2,\text{ for all }(t,\omega,y_1,y_2,z,a)$$
\item[\rm{(vii)}] $F$ is continuous in $a$.
\end{itemize}
\end{Assumption}

\vspace{0.8em}
\begin{Remark}\label{remark.debut}
Let us comment on the above assumptions. Assumptions \ref{assump.h} $\rm{(i)}$ and $\rm{(iv)}$ are taken from \cite{stz} and are needed to deal with the technicalities induced by the quasi-sure framework. Assumptions \ref{assump.h} $\rm{(ii)}$ and $\rm{(iii)}$ are quite standard in the classical BSDE literature. Then, Assumptions \ref{assump.h} $\rm{(v)}$ and $\rm{(vi)}$ were introduced by Pardoux in \cite{pard} in a more general setting (namely with a general growth condition in $y$, and only a continuity assumption on $y$) and are also quite common in the literature. Let us immediately point out that as explained in Remark \ref{remark.imp} below, we must restrict ourselves to linear growth in $y$, because of the technical difficulties due to the $2$BSDE framework. Moreover, we need to assume uniform continuity in $y$ to ensure that we have a strong convergence result for the approximation we will consider (see also Remark \ref{remark.imp}). Finally, Assumption \ref{assump.h} $\rm{(vii)}$ is needed in our framework to obtain technical results concerning monotone convergence in a quasi-sure setting.
\end{Remark}

\subsection{The spaces and norms}

We now recall from \cite{stz} the spaces and norms which will be needed for the formulation of the second order BSDEs. Notice that all subsequent notations extend to the case $\kappa=1$.

\vspace{0.8em}
For $p\geq 1$, $L^{p,\kappa}_H$ denotes the space of all $\mathcal F_T$-measurable scalar r.v. $\xi$ with
$$\No{\xi}_{L^{p,\kappa}_H}^p:=\underset{\mathbb{P} \in \mathcal{P}_H^\kappa}{\sup}\mathbb E^{\mathbb P}\left[|\xi|^p\right]<+\infty.$$

$\mathbb H^{p,\kappa}_H$ denotes the space of all $\mathbb F^+$-progressively measurable $\mathbb R^d$-valued processes $Z$ with
$$\No{Z}_{\mathbb H^{p,\kappa}_H}^p:=\underset{\mathbb{P} \in \mathcal{P}_H^\kappa}{\sup}\mathbb E^{\mathbb P}\left[\left(\int_0^T|\widehat a_t^{1/2}Z_t|^2dt\right)^{\frac p2}\right]<+\infty.$$

$\mathbb D^{p,\kappa}_H$ denotes the space of all $\mathbb F^+$-progressively measurable $\mathbb R$-valued processes $Y$ with
$$\mathcal P^\kappa_H-q.s. \text{ c\`adl\`ag paths, and }\No{Y}_{\mathbb D^{p,\kappa}_H}^p:=\underset{\mathbb{P} \in \mathcal{P}_H^\kappa}{\sup}\mathbb E^{\mathbb P}\left[\underset{0\leq t\leq T}{\sup}|Y_t|^p\right]<+\infty.$$

For each $\xi \in L^{1,\kappa}_H$, $\mathbb P\in \mathcal P^\kappa_H$ and $t \in [0,T]$ denote
$$\mathbb E_t^{H,\mathbb P}[\xi]:=\underset{\mathbb P^{'}\in \mathcal P^\kappa_H(t^{+},\mathbb P)}{\esup^\mathbb P}\mathbb E^{\mathbb P^{'}}_t[\xi] \text{ where } \mathcal P^\kappa_H(t^{+},\mathbb P):=\left\{\mathbb P^{'}\in\mathcal P^\kappa_H:\mathbb P^{'}=\mathbb P \text{ on }\mathcal F_t^+\right\}.$$

Here $\mathbb{E}_t^{\mathbb P}[\xi]:=E^{\mathbb P}[\xi|\mathcal F_t]$. Then we define for each $p\geq \kappa$,
$$\mathbb L_H^{p,\kappa}:=\left\{\xi\in L^{p,\kappa}_H:\No{\xi}_{\mathbb L_H^{p,\kappa}}<+\infty\right\} \text{ where } \No{\xi}_{\mathbb L_H^{p,\kappa}}^p:=\underset{\mathbb P\in\mathcal P^\kappa_H}{\sup}\mathbb E^{\mathbb P}\left[\underset{0\leq t\leq T}{\esup}^{\mathbb P}\left(\mathbb E_t^{H,\mathbb P}[|\xi|^\kappa]\right)^{\frac{p}{\kappa}}\right].$$

Finally, we denote by $\mbox{UC}_b(\Omega)$ the collection of all bounded and uniformly continuous maps $\xi:\Omega\rightarrow \mathbb R$ with respect to the $\No{\cdot}_{\infty}$-norm, and we let

$$\mathcal L^{p,\kappa}_H:=\text{the closure of $\mbox{UC}_b(\Omega)$ under the norm $\No{\cdot}_{\mathbb L^{p,\kappa}_H}$, for every $1\leq\kappa \leq p$}.$$

\subsection{Formulation}\label{sec.definition}
We shall consider the following second order BSDE (2BSDE for short), which was first defined in \cite{stz}

\begin{equation}
Y_t=\xi +\int_t^T\widehat{F}_s(Y_s,Z_s)ds -\int_t^TZ_sdB_s + K_T-K_t, \text{ } 0\leq t\leq T, \text{  } \mathcal P_H^\kappa-q.s.
\label{2bsde}
\end{equation}

\begin{Definition}
For $\xi \in L^{2,\kappa}_H$, we say $(Y,Z)\in \mathbb D^{2,\kappa}_H\times\mathbb H^{2,\kappa}_H$ is a solution to 2BSDE \reff{2bsde} if 

\begin{itemize}
\item[$\bullet$] $Y_T=\xi$, $\mathcal P_H^\kappa-q.s$. 
\item[$\bullet$] $\forall \mathbb P \in \mathcal P_H^\kappa$, the process $K^{\mathbb P}$ defined below has non-decreasing paths $\mathbb P-a.s.$ 
\begin{equation}
K_t^{\mathbb P}:=Y_0-Y_t - \int_0^t\widehat{F}_s(Y_s,Z_s)ds+\int_0^tZ_sdB_s, \text{ } 0\leq t\leq T, \text{  } \mathbb P-a.s.
\label{2bsde.k}
\end{equation}

\item[$\bullet$] The family $\left\{K^{\mathbb P}, \mathbb P \in \mathcal P_H^\kappa\right\}$ satisfies the minimum condition

\begin{equation}
K_t^{\mathbb P}=\underset{ \mathbb{P}^{'} \in \mathcal{P}_H(t^+,\mathbb{P}) }{ \einf^{\mathbb P} }\mathbb{E}_t^{\mathbb P^{'}}\left[K_T^{\mathbb{P}^{'}}\right], \text{ } 0\leq t\leq T, \text{  } \mathbb P-a.s., \text{ } \forall \mathbb P \in \mathcal P_H^\kappa.
\label{2bsde.min}
\end{equation}
\end{itemize}
Moreover if the family $\left\{K^{\mathbb P}, \mathbb P \in \mathcal P_H^\kappa\right\}$ can be aggregated into a universal process $K$, we call $(Y,Z,K)$ a solution of 2BSDE \reff{2bsde}.
\end{Definition}

\vspace{0.8em}
\begin{Remark}\label{rem.int}
Let us comment on this definition. As already explained, the PDE intuition leads us to think that the solution of a 2BSDE should be a supremum of solutions of standard BSDEs. Therefore, for each $\mathbb P$, the role of the non-decreasing process $K^\mathbb P$ is to "push" the process $Y$ so that it remains above the solution of the BSDE with terminal condition $\xi$ and generator $\widehat F$ under $\mathbb P$. In this regard, 2BSDEs share some similarities with reflected BSDEs. 

\vspace{0.8em}
Pursuing this analogy, the minimum condition \reff{2bsde.min} tells us that the processes $K^\mathbb P$ act in a "minimal" way (exactly as implied by the Skorohod condition for reflected BSDEs), and we will see in the next Section that it implies uniqueness of the solution. Besides, if the set $\mathcal P^\kappa_H$ was reduced to a singleton
$\{\mathbb P\}$, then \reff{2bsde.min} would imply that $K^\mathbb P$ is a martingale and a non-decreasing process and is therefore null. Thus we recover the standard BSDE theory.

\vspace{0.8em}
Finally, we would like to emphasize that in the language of G-expectation of Peng \cite{peng}, (\ref{2bsde.min}) is equivalent, at least if the family can be aggregated into a process $K$, to saying that $-K$ is a G-martingale. This has been already observed in \cite{stz3} where the authors proved the G-martingale representation property, which formally corresponds to a 2BSDE with a generator equal to $0$.
\end{Remark}

\vspace{0.8em}  
Following \cite{stz}, in addition to Assumption \ref{assump.h}, we will always assume

\begin{Assumption}\label{assump.h2}
\begin{itemize}
\item[\rm{(i)}] $\mathcal P_H^\kappa$ is not empty.
\item[\rm{(ii)}] The process $\widehat F^0 $ satisfies the following integrability condition
\begin{align}
\phi^{2,\kappa}_H&:=\underset{\mathbb P\in\mathcal P^\kappa_H}{\sup}\mathbb E^{\mathbb P}\left[\underset{0\leq t\leq T}{\esup}^{\mathbb P}\left(\mathbb E_t^{H,\mathbb P}\left[\int^T_0|\Fh^0_s|^\kappa ds\right]\right)^{\frac{2}{\kappa}}\right]<+\infty. 
\end{align}
\end{itemize}
\end{Assumption}

\vspace{0.8em}
Before going on, let us recall one of the main results of \cite{stz}. For this, we first recall their assumptions on the generator $F$
\begin{Assumption} \label{assump.stz}
\begin{itemize}
\item[\rm{(i)}] The domain $D_{F_t(y,z)}=D_{F_t}$ is independent of $(\omega,y,z)$.
\item[\rm{(ii)}] For fixed $(y,z,a)$, $F$ is $\mathbb{F}$-progressively measurable.
\item[\rm{(iii)}] We have the following uniform Lipschitz-type property
$$\forall (y,y',z,z',a,t,\omega), \text{ } \abs{F_t(\omega,y,z,a)- F_t(\omega,y',z',a)}\leq C\left(\abs{y-y'}+\abs{a^{1/2}(z-z')}\right).$$
\item[\rm{(iv)}] $F$ is uniformly continuous in $\omega$ for the $||\cdot||_\infty$ norm.
\end{itemize}
\end{Assumption}

\vspace{0.8em}
\begin{Theorem}[Soner, Touzi, Zhang \cite{stz}]
Let Assumptions \ref{assump.h2} and \ref{assump.stz} hold. Then, for any $\xi\in\mathcal L^{2,\kappa}_H$, the 2BSDE \reff{2bsde} has a unique solution $(Y,Z)\in\mathbb D^{2,\kappa}_H\times\mathbb H^{2,\kappa}_H$.
\end{Theorem}

\vspace{0.8em}
We now state the main result proved in this paper
\begin{Theorem}
\label{th1}
Suppose there exists $\epsilon>0$ such that $\xi \in \mathcal L^{2,\kappa}_H\cap L^{2+\epsilon,\kappa}_H$ and 

\begin{equation}\label{epsilon}
\underset{\mathbb P\in\mathcal P^\kappa_H}{\Sup}\mathbb E^\mathbb P\left[\int_0^T\abs{\widehat F^0_t}^{2+\epsilon}dt\right]<+\infty.
\end{equation}

\vspace{0.8em}
Then, under Assumptions \ref{assump.h} and \ref{assump.h2}, there exists a unique solution $(Y,Z)\in \mathbb D^{2,\kappa}_H\times \mathbb H^{2,\kappa}_H$ of the 2BSDE \reff{2bsde}.
\end{Theorem}

\vspace{0.8em}
\begin{Remark}
Notice that the above Theorems both correspond to the one dimensional case (in the sense that the terminal condition belongs to $\mathbb R$). This is different from the standard BSDE literature (see \cite{elkaroui} and \cite{pard}) where existence and uniqueness of a solution can be obtained in arbitrary dimensions under Lipschitz type assumptions on the generator. This difference is mainly due to the fact that the comparison theorem (only valid in dimension $1$) is intensely used to obtain existence and uniqueness for 2BSDEs in \cite{stz}. This strategy of proof is (partly) due to the fact that the classical fixed-point argument for standard BSDEs no longer works in the 2BSDE framework. Indeed, as can be seen from the estimates of Theorem $4.5$ in \cite{stz}, trying to reproduce the classical fixed-point argument, would only lead to an application which is $1/2$-H\"older continuous, which is not sufficient to conclude. Extending these results to the multidimensional case is nonetheless an important and difficult problem.
\end{Remark}

\vspace{0.8em}
\subsection{Representation and uniqueness of the solution}

We follow once more Soner, Touzi and Zhang \cite{stz}. For any $\mathbb{P}\in\mathcal{P}^\kappa_H$, $\mathbb{F}$-stopping time $\tau$, and $\mathcal{F}_\tau$-measurable random variable $\xi \in \mathbb{L}^2(\mathbb P)$, let $(y^\mathbb{P},z^\mathbb{P}):=(y^\mathbb{P}(\tau,\xi),z^\mathbb{P}(\tau,\xi))$ denote the unique solution to the following standard BSDE (existence and uniqueness under our assumptions follow from Pardoux \cite{pard})

\begin{equation}\label{bsdeee}
y_t^\mathbb P=\xi+\int_t^\tau\widehat F_s(y_s^\mathbb P,z_s^\mathbb P)ds-\int_t^\tau z_s^\mathbb PdB_s, \text{ }0\leq t\leq \tau,\text{ } \mathbb P-a.s.
\end{equation}

We then have similarly as in Theorem $4.4$ of \cite{stz}.

\begin{Theorem}\label{unique}
Let Assumptions \ref{assump.h} and \ref{assump.h2} hold. Assume $\xi \in \mathbb{L}^{2,\kappa}_H$ and that $(Y,Z)\in \mathbb D^{2,\kappa}_H\times\mathbb H^{2,\kappa}_H$ is a solution of the 2BSDE \reff{2bsde}. Then, for any $\mathbb{P}\in\mathcal{P}^\kappa_H$ and $0\leq t_1< t_2\leq T$,
\begin{equation}
Y_{t_1}=\underset{\mathbb{P}^{'}\in\mathcal{P}^\kappa_H(t_1^+,\mathbb{P})}{\esup^\mathbb{P}}y_{t_1}^{\mathbb{P}^{'}}(t_2,Y_{t_2}), \text{ }\mathbb{P}-a.s.
\label{representation}
\end{equation}
Consequently, the 2BSDE \reff{2bsde} has at most one solution in $ \mathbb D^{2,\kappa}_H\times\mathbb H^{2,\kappa}_H$.
\end{Theorem}

\vspace{0.8em}
Notice that the representation formula \reff{representation} corresponds exactly to the PDE intuition we described in the Introduction and Remark \ref{rem.int}, and shows that the solution to a 2BSDE is indeed a supremum over a family of standard BSDEs.

\vspace{0.8em}
\proof
The proof is almost identical to the proof of Theorem $4.4$ in \cite{stz}, and we reproduce it here for the convenience of the reader. First, if \reff{representation} holds, then
$$Y_{t}=\underset{\mathbb{P}^{'}\in\mathcal{P}^\kappa_H(t^+,\mathbb{P})}{\esup^\mathbb{P}}y_{t}^{\mathbb{P}^{'}}(T,\xi), \text{ } t\in [0,T], \text{ }\mathbb{P}-a.s. \text{ for all }\mathbb P\in \mathcal{P}^\kappa_H,$$
and thus is unique. Then, since we have that $d\left<Y,B\right>_t=Z_tdB_t, \text{ } \mathcal{P}^\kappa_H-q.s.$,  $Z$ is also unique. We shall now prove \reff{representation}.

\begin{itemize}
\item[\rm{(i)}] Fix $0\leq t_1<t_2\leq T$ and $\mathbb P\in\mathcal P^\kappa_H$. For any $\mathbb P^{'}\in\mathcal P^\kappa_H(t_1^+,\mathbb P)$, we have
$$Y_{t} = Y_{t_2} + \int_{t}^{t_2} \widehat{F}_s(Y_s,Z_s)ds - \int_{t}^{t_2} Z_sdB_s + K_{t_2}^{ \mathbb P^{'}} - K_{t}^{ \mathbb P^{'} }, \text{ } 0\leq t\leq t_2, \text{ } \mathbb P^{'}-a.s.$$
and that $K^{\mathbb P^{'}}$ is nondecreasing, $\mathbb P^{'}-a.s.$ Then, we can apply a generalized comparison theorem proved by Lepeltier, Matoussi and Xu (see Theorem $4.1$ in \cite{lmx}) under $\mathbb P^{'}$ to obtain $Y_{t_1}\geq y_{t_1}^{\mathbb P^{'}}(t_2,Y_{t_2})$, $\mathbb P^{'}-a.s.$ Since $\mathbb P^{'}=\mathbb P$ on $\mathcal F_t^+$, we get $Y_{t_1}\geq y_{t_1}^{\mathbb P^{'}}(t_2,Y_{t_2})$, $\mathbb P-a.s.$ and thus
$$Y_{t_1}\geq\underset{\mathbb{P}^{'}\in\mathcal{P}^\kappa_H(t_1^+,\mathbb{P})}{\esup^\mathbb{P}}y_{t_1}^{\mathbb{P}^{'}}(t_2,Y_{t_2}), \text{ }\mathbb{P}-a.s.$$

\item[\rm{(ii)}] We now prove the reverse inequality. Fix $\mathbb P\in\mathcal P^\kappa_H$. We will show in $\rm{(iii)}$ below that
$$C_{t_1}^\mathbb P:=\underset{\mathbb{P}^{'}\in\mathcal{P}^\kappa_H(t_1^+,\mathbb{P})}{\esup^\mathbb{P}}\mathbb E_{t_1}^{\mathbb{P}^{'}}\left[\left(K_{t_2}^{\mathbb P^{'}}-K_{t_1}^{\mathbb P^{'}}\right)^2\right]<+\infty,\text{ }\mathbb P-a.s.$$

For every $\mathbb P^{'}\in \mathcal P^\kappa_H(t^+,\mathbb P)$, denote
$$\delta Y:=Y-y^{\mathbb P^{'}}(t_2,Y_{t_2})\text{ and }\delta Z:=Z-z^{\mathbb P^{'}}(t_2,Y_{t_2}).$$

By the Lipschitz Assumption \ref{assump.h}$\rm{(iii)}$ and the monotonicity Assumption \ref{assump.h}$\rm{(vi)}$, there exist a bounded process $\lambda$ and a process $\eta$ which is bounded from above such that
$$\delta Y_t=\int_t^{t_2}\left(\eta_s\delta Y_s+\lambda_s\widehat{a}_s^{1/2}\delta Z_s\right)ds-\int_t^{t_2}\delta Z_sdB_s+K_{t_2}^{\mathbb P^{'}}-K_t^{\mathbb P^{'}}, \text{ }t\leq t_2,\text{ }\mathbb P^{'}-a.s.$$

Define for $t_1\leq t\leq t_2$
$$M_t:=\exp\left(\int_{t_1}^t\left(\eta_s-\frac12\abs{\lambda_s}^2\right)ds-\int_{t_1}^t\lambda_s\widehat{a}_s^{-1/2}dB_s\right),\text{ }\mathbb P^{'}-a.s.$$

By It\^o's formula, we obtain, as in \cite{stz}, that

$$\delta Y_{t_1}=\mathbb E_{t_1}^{\mathbb P^{'}}\left[\int_{t_1}^{t_2}M_tdK_t^{\mathbb P^{'}}\right]\leq\mathbb E_{t_1}^{\mathbb P^{'}}\left[\underset{t_1\leq t\leq t_2}{\sup}(M_t)(K_{t_2}^{\mathbb P^{'}}-K_{t_1}^{\mathbb P^{'}})\right],$$
since $K^{\mathbb P^{'}}$ is non-decreasing. Then, because $\lambda$ is bounded and $\eta$ is bounded from above, we have for every $p\geq 1$
$$\mathbb E_{t_1}^{\mathbb P^{'}}\left[\underset{t_1\leq t\leq t_2}{\sup}(M_t)^p\right]\leq C_p,\text{ } \mathbb P^{'}-a.s.$$
Then it follows from the H\"older inequality that
$$\delta Y_{t_1}\leq C(C_{t_1}^{\mathbb P^{'}})^{1/3}\left(\mathbb E_{t_1}^{\mathbb P^{'}}\left[K_{t_2}^{\mathbb P^{'}}-K_{t_1}^{\mathbb P^{'}}\right]\right)^{1/3},\text{ } \mathbb P^{'}-a.s.$$

By the minimum condition \reff{2bsde.min} and since $\mathbb P^{'}\in \mathcal P^\kappa_H(t^+,\mathbb P)$ is arbitrary, this ends the proof.

\item[\rm{(iii)}] It remains to show that the estimate for $C_{t_1}^{\mathbb P^{'}}$ holds. By definition of the family $\left\{K^\mathbb P, \mathbb P\in\mathcal P^\kappa_H\right\}$, the linear growth condition in $y$ of the generator and the Lipschitz condition in $z$, we have

$$\underset{\mathbb{P}^{'}\in\mathcal{P}^\kappa_H(t_1^+,\mathbb{P})}{\sup{^\mathbb{P}}}\mathbb E^{\mathbb{P}^{'}}\left[\left(K_{t_2}^{\mathbb P^{'}}-K_{t_1}^{\mathbb P^{'}}\right)^2\right]\leq C\left(1+\No{Y}^2_{\mathbb D^{2,\kappa}_H}+\No{Z}^2_{\mathbb H^{2,\kappa}_H}+\phi^{2,\kappa}_H\right)<+\infty.$$

Then we proceed exactly as in the proof of Theorem $4.4$ in \cite{stz}. 
\end{itemize}
\ep

\subsection{Non-dominated family of probability measures and monotone Convergence Theorem}

Let $\mathcal P_0\subset\overline{\mathcal P}_W$. In general in a non-dominated framework, the monotone convergence Theorem may not hold, that is to say that even if we have a sequence of random variables $X_n$ which decreasingly converges $\mathcal P_0-q.s.$ to $0$, then we may not have that
$$\underset{\mathbb P\in\mathcal P_0}{\sup}\mathbb E^\mathbb P[X_n]\downarrow0.$$

\vspace{0.8em}
Indeed, let us consider for instance the set 
$$\mathcal P_1:=\left\{\mathbb P^p:=\mathbb P_0\circ (\sqrt{p}B),\text{ }p\in\mathbb N^*\right\}.$$

\vspace{0.8em}
Then, define $Y_n:=\frac{B_1^2}{n}$. It is clear that the sequence $Y_n$ decreases $\mathbb P$-a.s. to $0$ for all $\mathbb P\in\mathcal P_1$. But we have for all $p$ and all $n$
$$\mathbb E^{\mathbb P^p}\left[Y_n\right]=\frac{p}{n},$$
which implies that $\underset{\mathbb P\in\mathcal P_1}{\sup}\mathbb E^\mathbb P[Y_n]=+\infty$.

\vspace{0.8em}
It is therefore clear that it is necessary to add more assumptions in order to recover a monotone convergence theorem. Notice that those assumptions will concern both the family $\mathcal P_0$ and the random variables considered. For instance, in the above example, this is the fact that the set $\mathcal P_1$ considered is not weakly compact which implies that the monotone convergence theorem can fail.

\vspace{0.8em}
Let us now provide some definitions

\begin{Definition}
Let $(X^\mathbb P)_{\mathbb P\in\mathcal P_0}$ be a family of random variables. The family can be aggregated if there exists a random variable $X$ such that
$$X=X^\mathbb P,\ \mathbb P-a.s.,\ \text{for all }\mathbb P\in\mathcal P_0.$$

\vspace{0.8em}
$X$ will be called an aggregator for this family.
\end{Definition}

\vspace{0.8em}
\begin{Definition}
\begin{itemize}
\item[$\rm{(i)}$] We say that a family of random variables $(X^\mathbb P)_{\mathbb P\in\mathcal P_0}$ is $\mathcal P_0$-uniformly integrable if
$$\underset{C\rightarrow+\infty}{\lim}\text{ }\underset{\mathbb P\in\mathcal P_0}{\sup}\mathbb E^\mathbb P\left[\abs{X^\mathbb P}1_{\abs{X^\mathbb P}>C}\right]=0.$$
\item[$\rm{(ii)}$] We say that a family of random variables $(X^\mathbb P)_{\mathbb P\in\mathcal P_0}$ is $\mathcal P_0$-quasi continuous if for all $\mathbb P\in\mathcal P_0$ and for all $\epsilon>0$, there exists an open set $\mathcal O^{\epsilon}\subset\Omega$ such that
$$\mathbb P(\mathcal O^{\epsilon})\leq \epsilon\text{ and $X^\mathbb P$ is continuous in $\omega$ outside $\mathcal O^{\epsilon}$}.$$
\end{itemize}
\end{Definition}

\vspace{0.8em}
We will keep the same terminology if the family $(X^\mathbb P)_{\mathbb P\in\mathcal P_0}$ can be aggregated.

\vspace{0.8em}
\begin{Remark}
\begin{itemize}
\item A sufficient condition for a family of random variables $(X^\mathbb P)_{\mathbb P\in\mathcal P_0}$ to be $\mathcal P_0$-uniformly integrable is that 
$$\underset{\mathbb P\in\mathcal P_0}{\sup}\mathbb E^\mathbb P\left[\abs{X^\mathbb P}^{1+\epsilon}\right]<+\infty,\text{ for some }\eps>0.$$
\item The definition of $\mathcal P_0$-quasi continuity is inspired by the notion of quasi-continuity given in \cite{dhp}, but adapted to our context where we work without the theory of capacities and where the random variables are not necessarily aggregated.
\end{itemize}
\end{Remark}

\vspace{0.8em}
We then have the following monotone convergence theorem proved in \cite{dhp}

\begin{Theorem}[Denis, Hu, Peng]\label{th.monopeng}
Let $\mathcal P_0$ be a weakly compact family and let $X_n$ be a sequence of $\mathcal P_0$-quasi continuous and $\mathcal P_0$-uniformly integrable random variables which verifies that
$$X_n(\omega)\downarrow 0, \ \text{for every $\omega\in\Omega\backslash\mathcal N$, where $\underset{\mathbb P\in\mathcal P_0}{\sup}\mathbb P(\mathcal N)=0$}.$$

\vspace{0.8em}
Then
$$\underset{\mathbb P\in\mathcal P_0}{\sup}\mathbb E^\mathbb P[X_n]\downarrow0.$$
\end{Theorem}

\vspace{0.8em}
It is obvious that this Theorem is particularly suited for a capacity framework, as considered in \cite{denis} or \cite{dhp}. However in our case, we do not work with capacities, and all our properties hold only $\mathbb P$-a.s. for all probability measures in $\mathcal P_0$, which generally makes this Theorem not general enough for our purpose. Nonetheless, we will see later that Theorem \ref{th.monopeng} will only be applied to quantities which are not only defined for every $\omega$, but are also continuous in $\omega$. In that case, this difference between capacities and probabilities is not important.

\section{Proof of existence}\label{section.2cont}
\subsection{Preliminary results}
In order to prove existence, we need an approximation of continuous functions by Lipschitz functions proved by Lepeltier and San Martin in \cite{lsm}. %In the present framework, we also use the well-known approximation of Moreau and Yosida (see \cite{m} and \cite{y})

\vspace{0.8em}
Define $$ F_t^n(y,z,a):=\underset{u\in\mathbb Q}{\Inf}\left\{ F_t(u,z,a)+n\abs{y-u}\right\},$$
where the minimum over rationals $u\in\mathbb Q$ coincides with the minimum over real parameters by our continuity assumptions, and ensures the measurability of $F^n$.

\begin{Lemma} \label{lemma.approx}
Let $C$ be the constant in Assumption \ref{assump.h}. We have
\begin{itemize}
\item[\rm{(i)}] $F^n$ is well defined for $n\geq C$ and we have
$$\abs{F_t^n(y,z,a)}\leq \abs{ F_t(0,0,a)}+ C(1+\abs{y}+ |a^{1/2}z|), \text{ for all }(y,z,a,n,t,\omega).$$
\item[\rm{(ii)}] $|F_t^n(y,z_1,a)-F_t^n(y,z_2,a)|\leq C|a^{1/2}(z_1-z_2)|,$ for all $(y,z_1,z_2,a,t,\omega)$.
\item[\rm{(iii)}] $\abs{F_t^n(y_1,z,a)-F_t^n(y_2,z,a)}\leq n\abs{y_1-y_2}$, for all $(y_1,y_2,z,a,t,\omega)$.
\item[\rm{(iv)}] The sequence $(F_t^n(y,z,a))_{n}$ is increasing, for all $(t,y,z,a)$.
\item[\rm{(v)}] If $F$ is decreasing in $y$, then so is $F^n$.
\end{itemize}
\end{Lemma}

\vspace{0.8em}
\proof
The properties $\rm{(i)}$, $\rm{(ii)}$, $\rm{(iii)}$ and $\rm{(iv)}$ are either clear or proved in \cite{lsm}, thus we concentrate on $\rm{(v)}$.

\vspace{0.8em}
We assume that $F$ is decreasing in $y$. In particular, $F$ is differentiable in $y$ for $a.e.$ $y$. Define for all $y$, $h_{n,y,z,t,a}(u):= F_t(u,z,a)+n|y-u|$. For $u\leq y$, $h_{n,y,z,t,a}$ is clearly decreasing in $u$. Then, its minimum in $u$ can only be attained at $y$ or at a point strictly greater than $y$. 

\vspace{0.8em}
Therefore we can write

\begin{align*}
F_t^n(y,z,a)&=\min\left\{ F_t(y,z,a), \underset{u\in \mathbb Q,\text{ }u>y}{\Inf}\left\{ F_t(u,z)+n(u-y)\right\}\right\}\\
&=\min\left\{F_t(y,z,a), \underset{u\in \mathbb Q,\text{ }u>0}{\Inf}\left\{ F_t(u+y,z,a)+nu\right\}\right\},
\end{align*}
and under this form it is clear that $F^n$ is decreasing in $y$.
\ep

%\vspace{0.8em}
%\begin{Remark}
%In Theorem $5.1$ in Chapter $1$ of \cite{clark}, the function considered is supposed to be bounded from below. However, a careful reading of the proof shows that this hypothesis is only necessary to prove that the inf-convolution is well defined, which we already know to hold true is our case for $n$ large enough, thanks to our linear growth hypothesis.
%\end{Remark}
%
%\vspace{0.8em}
%\begin{Remark}
%Unlike \cite{lsm} or \cite{mat}, we do not use linear inf-convolution for our approximation, but a mix of linear and quadratic inf-convolution. This is due to our crucial need that our approximation remains uniformly Lipschitz in $z$ in the sense of Assumption \ref{assump.h}$\rm{(iii)}$. It is then more convenient to use quadratic inf-convolution in the variable $z$ to be able to use the results of non-smooth analysis of \cite{clark}.
%\end{Remark}

\vspace{0.8em}
We note that in above lemmas, and in all subsequent results, we shall denote by $C$ a generic constant which may vary from line to line and depends only on the dimension $d$, the maturity $T$ and the constants in Assumptions \ref{assump.h} and \ref{assump.h2}. We shall also denote by $C_{\kappa}$ a constant which may depend on $\kappa$ as well.

\vspace{0.8em}
Let us now note that we can always consider without loss of generality that the constant $\mu$ in Assumption \ref{assump.h}$\rm{(vi)}$ is equal to $0$. Indeed, we have the following Lemma

\begin{Lemma}\label{choubidou}
Let $\lambda>0$, then $(Y_t,Z_t,\left\{K_t^\mathbb P,\mathbb P\in \mathcal P^\kappa_H\right\})$ solve the 2BSDE \reff{2bsde} if and only if $\left(e^{\lambda t}Y_t,e^{\lambda t}Z_t,\left\{\int_0^t e^{\lambda s}dK_s^{\mathbb P},\mathbb P\in\mathcal P_H^\kappa\right\}\right)$ solve the 2BSDE with terminal condition $\overline{\xi}:=e^{\lambda T}\xi$ and driver $\widetilde F^{(\lambda)}_t(y,z):= F^{(\lambda)}_t(y,z,\widehat a_t)$, where $ F^{(\lambda)}_t(y,z,a):=e^{\lambda t} F_t (e^{-\lambda t} y, e^{-\lambda t} z,a)-\lambda y$.
\end{Lemma}

\vspace{0.8em}
\proof
The fact that the two solutions solve the corresponding equations is a simple consequence of It\^o's formula. The only thing that we have to check is that the family $\left\{K_t^\mathbb P,\mathbb P\in \mathcal P^\kappa_H\right\}$ satisfies the minimum condition \reff{2bsde.min} if and only if it is verified by the family $\left\{\int_0^t e^{\lambda s}dK_s^{\mathbb P},\mathbb P\in\mathcal P_H^\kappa\right\}$. 

\vspace{0.8em}
First of all, it is clear that 
$$\int_0^te^{\lambda s}dK_s^{\mathbb P}=\underset{\mathbb P^{'}\in\mathcal P^\kappa_H(t^+,\mathbb P)}{\einf^\mathbb P}\mathbb E_t^{\mathbb P^{'}}\left[\int_0^Te^{\lambda s}dK_s^{\mathbb P^{'}}\right]\Leftrightarrow \underset{\mathbb P^{'}\in\mathcal P^\kappa_H(t^+,\mathbb P)}{\einf^\mathbb P}\mathbb E_t^{\mathbb P^{'}}\left[\int_t^Te^{\lambda s}dK_s^{\mathbb P^{'}}\right]=0.$$

\vspace{0.8em}
Now for every $t\in [0,T]$, $\mathbb P\in\mathcal P^\kappa_H$ and $\mathbb P^{'}\in\mathcal P^\kappa_H(t^+,\mathbb P)$, the result follows from 
\begin{align*}
&\mathbb E_t^{\mathbb P^{'}}\left[\int_t^Te^{\lambda (s-T)}dK_s^{\mathbb P^{'}}\right]\leq \mathbb E_t^{\mathbb P^{'}} \left[K_T^{\mathbb P^{'}}-K_t^{\mathbb P}\right]\leq \mathbb E_t^{\mathbb P^{'}}\left[\int_t^Te^{\lambda s}dK_s^{\mathbb P^{'}}\right],\text{ }\mathbb P-a.s.
\end{align*}
\ep

\vspace{0.8em}
Thus, if we choose $\lambda=\mu$ then $F^{(\mu)}$ satisfies
$$(y_1-y_2)(F^{(\mu)}_t(\omega,y_1,z,a)-F^{(\mu)}_t(\omega,y_2,z,a))\leq 0,\text{ for all }(t,\omega,y_1,y_2,z,a).$$

\vspace{0.8em}
As a consequence of Lemma \ref{choubidou}, we can assume without loss of generality that our driver is decreasing in $y$. Therefore, frow now on this assumption will replace Assumption \ref{assump.h}$\rm{(vi)}$. 

\vspace{0.8em}
As explained in Remark \ref{remark.imp}, we will actually need a strong convergence result for the sequence $$\widehat F^n_t(y,z):=F^n_t(y,z,\widehat a_t).$$

\vspace{0.8em}
Let us define the following quantity
\begin{align*}
&\widetilde  F^n_t:=\underset{(y,z,a)\in\mathbb R^{d+1}\times[\underline{a},\bar{a}]}{\sup}\left\{F_t(y,z,a)-F_t^n(y,z,a)\right\}.
\end{align*}

\vspace{0.8em}
We then have the following result

\vspace{0.8em}
\begin{Lemma}
\label{assump.conv}
Let Assumption \ref{assump.h} hold. Then the sequence $\widehat F^n$ converges uniformly globally in $(y,z)$ and for all $0\leq t\leq T$ and all $\eps>0$
$$\underset{\mathbb P\in\mathcal P^\kappa_H}{\sup}\mathbb E^\mathbb P\left[\abs{\widetilde  F^{n}_t}^{2+\epsilon}\right]=\underset{\mathbb P\in\overline{\mathcal P}^\kappa_H}{\sup}\mathbb E^\mathbb P\left[\abs{\widetilde  F^{n}_t}^{2+\epsilon}\right]\leq C,$$
 for some $C$ independent of $n$.
\end{Lemma}

\vspace{0.8em}
The proof is relegated to the Appendix.

\vspace{0.8em}
\begin{Remark}
Notice that in the above Lemma, we emphasize that we consider a supremum over $\overline{\mathcal P}^\kappa_H$ (even though it is equal to the supremum over $\mathcal P^\kappa_H$). This is important because we are going to apply the monotone convergence Theorem of \cite{dhp} to the quantity $\abs{\widetilde  F^{n}_t}$ in the sequel, and this theorem can only be used under a weakly compact family of probability measures.
\end{Remark}

\vspace{0.8em}
We are now in a position to prove the main result of this paper Theorem \ref{th1}

\subsection{Proof of the main result}
For a fixed $n$, consider the following 2BSDE
\begin{equation}
Y_t^n=\xi +\int_t^T\widehat{F}_s^n(Y_s^n,Z_s^n)ds -\int_t^TZ_s^ndB_s + K_T^n-K_t^n, \text{ } 0\leq t\leq T, \text{  } \mathcal P_H^\kappa-q.s.
\label{2bsden}
\end{equation}

\vspace{0.8em}
By Lemma \ref{lemma.approx} and our Assumptions \ref{assump.h} and \ref{assump.h2} we know that all the requirements of Theorem $4.7$ of \cite{stz} are fulfilled. Thus, we know that for all $n$ the above 2BSDE has a unique solution $(Y^n,Z^n)\in \mathbb D^{2,\kappa}_H\times \mathbb H^{2,\kappa}_H$. Moreover, if we introduce the following standard BSDEs for all $\mathbb P \in \mathcal P_H^\kappa$
\begin{equation}
y_t^{\mathbb P,n}=\xi +\int_t^T\widehat{F}_s^n(y_s^{\mathbb P,n},z_s^{\mathbb P,n})ds -\int_t^Tz_s^{\mathbb P,n}dB_s, \text{ } 0\leq t\leq T, \text{  } \mathbb P-a.s.,
\label{bsden}
\end{equation}
we have the already mentioned representation (see Theorem $4.4$ in \cite{stz})
\begin{equation}
Y_t^{n}=\underset{\mathbb{P}^{'}\in \mathcal{P}_H^\kappa(t^+,\mathbb{P})}{\esup^{\mathbb{P}}}y_t^{\mathbb{P}^{'},n}, \text{ } 0\leq t\leq T, \text{  } \mathbb{P}-a.s., \text{  } \forall \mathbb{P} \in \mathcal{P}^\kappa_H.
\label{rep}
\end{equation}

\vspace{0.8em}
The idea of the proof of existence is to prove that the limit in a certain sense of the sequence $(Y^n,Z^n)$ is a solution of the 2BSDE \reff{2bsde}. We first provide a priori estimates which are uniform in $n$ on the solutions of \reff{2bsden} and \reff{bsden}.

\vspace{2em}
\begin{Lemma}\label{estim}
There exists a constant $C_\kappa>0$ such that for all $n$ large enough
\begin{align*}
\No{Y^n}_{\mathbb D^{2,\kappa}_H}^2 + \No{Z^n}_{\mathbb H^{2,\kappa}_H}^2 + \underset{\mathbb P \in \mathcal P_H^\kappa}{\Sup}\mathbb{E}^{\mathbb P}\left[\abs{K^{\mathbb P,n}_T}^2+\underset{0\leq t\leq T}{\Sup}|y_t^{\mathbb P,n}|^2+\int_0^T|\widehat{a}_s^{1/2}z_s^{\mathbb P,n}|^2ds\right]\\
\leq C_{\kappa}\left(1+\No{\xi}^2_{\mathcal L^{2,\kappa}_H}+\phi^{2,\kappa}_H\right).
\end{align*}
\end{Lemma} 

\vspace{0.8em}
\proof
Let us consider the following BSDE
\begin{equation}
u_t^{\mathbb P}=\abs{\xi} +\int_t^T |\widehat F^0_s|+C\left(1+|u_s^{\mathbb P}|+|\widehat{a}_s^{1/2}v_s^{\mathbb P}|\right)ds -\int_t^Tv_s^{\mathbb P}dB_s, \text{ } 0\leq t\leq T, \text{  } \mathbb P-a.s.
\label{bsde}
\end{equation}

\vspace{0.8em}
Since its generator is clearly Lipschitz it has a unique solution. Moreover, we can apply the comparison theorem of El Karoui, Peng and Quenez \cite{elkaroui} to obtain that, due to our uniform growth assumption and $\rm{(i)},\rm{(ii)},\rm{(iii)}$ of Lemma \ref{lemma.approx} 
$$\forall \text{ } m\leq n \text{ large enough, }\forall \text{ } \mathbb P\in \mathcal P_H^\kappa \text{, }  y^{\mathbb P,m}\leq y^{\mathbb P,n}\leq u^{\mathbb P}, \text{ } \mathbb P-a.s.$$

\vspace{0.8em}
Now, following line-by-line the proof of Lemma $4.3$ in \cite{stz} we obtain that for all $\mathbb P\in \mathcal P_H^\kappa$ and all $0\leq t\leq T$
$$|u_t^{\mathbb P}|\leq C_\kappa\left(1+\mathbb{E}_t^{\mathbb P}\left[|\xi|^\kappa+\int_t^T|\widehat F^0_s|^\kappa ds\right]^{1/\kappa}\right).$$

\vspace{0.8em}
Therefore by definition of the norms and the representation \reff{rep} we have
$$\No{Y^n}_{\mathbb D^{2,\kappa}_H}^2 + \underset{\mathbb P \in \mathcal P_H^\kappa}{\Sup}\mathbb{E}^{\mathbb P}\left[\underset{0\leq t\leq T}{\Sup}|y_t^{\mathbb P,n}|^2\right]\leq C_{\kappa}\left(1+\No{\xi}^2_{\mathcal L^{2,\kappa}_H}+\phi^{2,\kappa}_H\right).$$

\vspace{0.8em}
We next apply It\^o's formula to $(Y_t^n)^2$ under each $\mathbb P \in \mathcal P_H^\kappa$, we have $\mathbb P-a.s.$
$$(Y_t^n)^2 +\int_t^T|\widehat{a}^{1/2}_sZ_s^n|^2ds\leq\xi^2+2\int_t^TY^n_s\widehat{F}_s^n(Y^n_s,Z_s^n)ds - 2\int_t^TY_{s^-}^nZ_s^ndB_s+2\int_t^TY_{s^-}^ndK_s^{\mathbb P,n}.$$

\vspace{0.8em}
Thus, since $(Y^n,Z^n)\in \mathbb D^{2,\kappa}_H\times \mathbb H^{2,\kappa}_H$, by taking expectation we obtain that for all  $\mathbb P \in \mathcal P_H^\kappa$,
$$\mathbb{E}^{\mathbb P}\left[\int_t^T|\widehat{a}^{1/2}_sZ_s^n|^2ds\right]\leq \No{\xi}^2_{\mathbb L^{2,\kappa}_H}+2\mathbb{E}^{\mathbb P}\left[\int_t^TY^n_s\widehat{F}_s^n(Y^n_s,Z_s^n)ds +\underset{0\leq t\leq T}{\Sup}|Y_t^n|K_T^{\mathbb P,n}\right].$$

\vspace{0.8em}
Now the uniform growth condition $\rm{(i)}$ of Lemma \ref{lemma.approx} and the elementary inequality $2ab\leq a^2/\eps + \eps b^2,$ $\forall$ $\eps>0$ yield
\begin{align*}
\mathbb{E}^{\mathbb{P}} \left[ \int_t^T |\widehat{a}^{1/2}_sZ_s^n|^2ds\right]& \leq \No{\xi}^2_{\mathbb{L}^{2,\kappa}_H}+C\left(1+\phi^{2,\kappa}_H+\mathbb{E}^{\mathbb{P}} \left[\int_t^T |Y^n_s|^2ds\right]\right)\\
&\hspace{0.9em}+\frac13\mathbb{E}^{\mathbb{P}}\left[ \int_t^T |\widehat{a}_s^{1/2}Z^n_s|^2ds\right]+2\mathbb{E}^{\mathbb P}\left[\underset{0\leq t\leq T}{\Sup} |Y_t^n| K_T^{\mathbb{P},n}\right]\\
i.e.\text{ }\frac23  \mathbb{E}^{\mathbb{P}} \left[ \int_0^T |\widehat{a}^{1/2}_sZ_s^n|^2ds\right] &\leq \No{\xi}^2_{\mathbb{L}^{2,\kappa}_H} + C\left(1+(1+\eps^{-1})\No{Y^n}^2_{\mathbb{D}^{2,\kappa}_H}\right) +\eps \mathbb{E}^{\mathbb{P}} \left[(K_T^{\mathbb{P},n})^2\right].
\end{align*}

\vspace{0.8em}
But by definition of $K_T^{\mathbb P,n}$, it is clear that
\begin{equation}\label{K}
\mathbb{E}^{\mathbb{P}} \left[(K_T^{\mathbb{P},n})^2\right]\leq C_0\left(1+\No{\xi}^2_{\mathbb{L}^{2,\kappa}_H}+\phi^{2,\kappa}_H+\No{Y^n}^2_{\mathbb{D}^{2,\kappa}_H} + \int_0^T |\widehat{a}^{1/2}_sZ_s^n|^2ds\right).
\end{equation}

\vspace{0.8em}
Choosing $\eps=\frac{1}{3C_0}$, reporting \reff{K} in the previous inequality and taking supremum over $\mathbb P$ yields
$$\No{Z^n}_{\mathbb H^{2,\kappa}_H}^2\leq C_{\kappa}\left(1+\No{\xi}^2_{\mathbb L^{2,\kappa}_H}+\phi^{2,\kappa}_H\right),$$
from which we can then deduce the result for $K_T^{\mathbb P,n}$.

\vspace{0.8em}
Finally, we can show similarly, by applying It\^o's formula to $y_t^{\mathbb P,n}$ instead, that
$$\underset{\mathbb P \in \mathcal P_H^\kappa}{\Sup}\mathbb{E}^{\mathbb P}\left[\int_0^T|\widehat{a}_s^{1/2}z_s^{\mathbb P,n}|^2ds\right]\leq C_{\kappa}\left(1+\No{\xi}^2_{\mathbb L^{2,\kappa}_H}+\phi^{2,\kappa}_H\right).$$
\ep

\vspace{0.8em}
Now from the comparison theorems for 2BSDEs (see Corolary $4.5$ in \cite{stz}) and BSDEs and $\rm{(iv)}$ of Lemma \ref{lemma.approx}, we recall that we have for all $0\leq t \leq T$
$$y^{\mathbb P,n}_t\leq y^{\mathbb P,n+1}_t \text{ }\mathbb P-a.s. \text{ for all } \mathbb P\in \mathcal P_H^\kappa \text{, and } Y^n_t\leq Y^{n+1}_t, \text{ } \mathcal P_H^\kappa-q.s.$$

\vspace{0.8em}
\begin{Remark}\label{remark.imp}
If we were in the classical framework, this $\mathbb P-$almost sure convergence of $y^{\mathbb P,n}$ together with the estimates of Lemma \ref{estim} would be sufficient to prove the convergence in the usual $\mathbb H^2$ space, thanks to the dominated convergence theorem. However, in our case, since the norms involve the supremum over a family of probability measures, this theorem can fail. This is exactly the major difficulty when considering the 2BSDE framework, since most of the techniques used in the standard BSDE literature to prove existence results involve approximations. In order to solve this problem, we need more regularity to be able apply the monotone convergence Theorem $31$ of \cite{dhp}. This is exactly why we had to add the assumption of uniform continuity in $y$ for our proof to work. Moreover, this also explains why, as already mentioned in Remark \ref{remark.debut}, we cannot generalize completely the results of Pardoux \cite{pard}. Indeed, restricting ourselves to linear growth in $y$ allows us to use the approximation by inf-convolution which has some very nice properties. If we had considered general growth in $y$, then it would have been extremely difficult to find reasonable conditions on the driver $\widehat F$ in order to have uniform convergence of the approximation.
\end{Remark}

\vspace{0.8em}
Next, we prove that

\begin{Lemma}
\label{lemma.yl2}
Under the hypotheses of Theorem \ref{th1}, we have
\begin{equation*}
\underset{0\leq t\leq T}{\Sup}\ \underset{\mathbb P \in \mathcal P_H^\kappa}{\Sup}\mathbb{E}^{\mathbb P}\left[|y_t^{\mathbb P}-y_t^{\mathbb P,n}|^{2+\epsilon'}\right]\underset{n\rightarrow +\infty}{\longrightarrow}0,
\end{equation*}
for any $0<\epsilon'<\epsilon$, with the same $\epsilon$ as in \reff{epsilon}.
\end{Lemma}

\proof
Thanks to Lemma \ref{lemma.approx}, we know that we can apply the same proof as that of Theorem $2.3$ of \cite{pard} in order to get for each $\mathbb P\in\mathcal P^\kappa_H$
\begin{align*}
\mathbb E^{\mathbb P}\left[|y_t^{\mathbb P}-y_t^{\mathbb P,n}|^{2+\epsilon'}\right]&\leq C\int_0^T\mathbb E^{\mathbb P}\left[|\widehat F_s(y_s^{\mathbb P},z_s^{\mathbb P})-\widehat F_s^n(y_s^{\mathbb P},z_s^{\mathbb P})|^{2+\epsilon'}\right]ds\leq C\int_0^T\mathbb E^{\mathbb P}\left[|\widetilde F_s^n|^{2+\epsilon'}\right]ds.
\end{align*}

\vspace{0.8em}
By Lemma \ref{assump.conv}, we know that $F^n$ converges uniformly in $(y,z)$. Since $F^n$ and $F$ are also continuous in $a$ by Assumption \ref{assump.h}$\rm{(vii)}$, the convergence is also uniform in $a\in[\underline a,\bar a]$ by Dini's lemma. Thus, $|\widetilde F_s^n|^{2+\epsilon'}$ decreases to $0$ for every $\omega\in\Omega$.

\vspace{0.8em}
Then, since $F$ and $F^n$ are uniformly continuous in $\omega$ on the whole space $\Omega$, then $\abs{\widetilde  F^{n}_t}^{2+\epsilon'}$ is continuous in $\omega$ on $\Omega$, and therefore quasi-continuous in the sense of \cite{dhp}. Moreover, we have again by Lemma \ref{assump.conv} and the fact that $\epsilon'<\epsilon$
$$\underset{\mathbb P\in\overline{\mathcal P}_H^\kappa}{\sup}\mathbb E^\mathbb P\left[\left(\abs{\widetilde  F^{n}_t}^{2+\epsilon'}\right)^{1+\epsilon''}\right],$$
for some $\epsilon''>0$.

\vspace{0.8em}
Hence, we have by classical arguments (namely H\"older and Markov inequalities) that
$$\underset{N\rightarrow+\infty}{\lim}\ \underset{\mathbb P\in\overline{\mathcal P}^\kappa_H}{\sup}\mathbb E^\mathbb P\left[\abs{\widetilde  F^{n}_t}^{2+\epsilon'}1_{\abs{\widetilde  F^{n}_t}^{2+\epsilon'}>N}\right]=0.$$

\vspace{0.8em}
Therefore, we can apply the monotone convergence Theorem \ref{th.monopeng} under the family $\overline{\mathcal P}^\kappa_H$ to obtain that 
$$\underset{\mathbb P\in\mathcal P^\kappa_H}{\sup}\mathbb E^\mathbb P\left[\abs{\widetilde  F^{n}_t}^{2+\epsilon'}\right]=\underset{\mathbb P\in\overline{\mathcal P}^\kappa_H}{\sup}\mathbb E^\mathbb P\left[\abs{\widetilde  F^{n}_t}^{2+\epsilon'}\right]\underset{n\rightarrow +\infty}{\longrightarrow}0,\text{ for all }0\leq t\leq T.$$

\vspace{0.8em}
Finally, the required result follows from the standard dominated convergence Theorem for the integral with respect to the Lebesgue measure.
\ep

\vspace{0.8em}
We continue with the following result 
\begin{Lemma}\label{conv.d2y}
Assume moreover that there exists an $\epsilon>0$ such that $\xi\in L^{2+\epsilon,\kappa}_H$ and 
$$\underset{\mathbb P\in\mathcal P^\kappa_H}{\Sup}\mathbb E^\mathbb P\left[\int_0^T\abs{\widehat F^0_t}^{2+\epsilon}dt\right]<+\infty.$$

Then, we have
\begin{equation*}
\underset{\mathbb P \in \mathcal P_H^\kappa}{\Sup}\mathbb{E}^{\mathbb P}\left[\underset{0\leq t \leq T}{\Sup}|Y_t^n-Y_t^p|^2 \right]\longrightarrow0,\text{ as $n,p\longrightarrow +\infty$.}
\end{equation*}
\end{Lemma}

\vspace{0.8em}
\proof
By the representation \reff{rep}, we have for all $n,p$ large enough
$$\underset{0\leq t \leq T}{\Sup}|Y_t^n-Y_t^p|\leq \underset{0\leq t \leq T}{\Sup}\ \underset{\mathbb{P}^{'}\in \mathcal{P}_H^\kappa(t^+,\mathbb{P})}{\esup^{\mathbb{P}}}|y_t^{\mathbb{P}^{'},n}-y_t^{\mathbb{P}^{'},p}|, \text{ } \mathbb P-a.s. \text{ for all } \mathbb P\in \mathcal P_H^\kappa.$$

\vspace{0.8em}
Then, we easily get
$$\underset{0\leq t \leq T}{\Sup}|Y_t^n-Y_t^p|^2\leq \underset{0\leq t \leq T}{\Sup}\left(\underset{\mathbb{P}^{'}\in \mathcal{P}_H^\kappa(t^+,\mathbb{P})}{\esup^{\mathbb{P}}}\mathbb E_t^{\mathbb P^{'}}\left[\underset{0\leq s\leq T}{\Sup}|y_s^{\mathbb{P}^{'},n}-y_s^{\mathbb{P}^{'},p}|\right]\right)^2, \text{ } \mathbb P-a.s. \text{ for all } \mathbb P\in \mathcal P_H^\kappa.$$

\vspace{0.8em}
Taking expectations yields for all $\mathbb P\in \mathcal P^\kappa_H$

$$\mathbb E^\mathbb P\left[\underset{0\leq t \leq T}{\Sup}|Y_t^n-Y_t^p|^2\right]\leq \underset{\mathbb P\in\mathcal P^\kappa_H}{\Sup}\mathbb E^\mathbb P\left[\underset{0\leq t \leq T}{\Sup}\left(\underset{\mathbb{P}^{'}\in \mathcal{P}_H^\kappa(t^+,\mathbb{P})}{\esup^{\mathbb{P}}}\mathbb E_t^{\mathbb P^{'}}\left[\underset{0\leq s\leq T}{\Sup}|y_s^{\mathbb{P}^{'},n}-y_s^{\mathbb{P}^{'},p}|\right]\right)^2\right].$$

\vspace{0.8em}
We next use the generalization of Doob maximal inequality of Proposition \ref{doob} in the Appendix (see also Lemma $6.2$ in \cite{stz3}), to obtain that for all $\epsilon'<\epsilon$
$$\mathbb E^\mathbb P\left[\underset{0\leq t \leq T}{\Sup}|Y_t^n-Y_t^p|^2\right]\leq C\underset{\mathbb P\in\mathcal P^\kappa_H}{\Sup}\mathbb E^\mathbb P\left[\underset{0\leq s\leq T}{\Sup}|y_s^{\mathbb{P},n}-y_s^{\mathbb{P},p}|^{2+\epsilon'}\right]^{\frac{2}{2+\epsilon'}}.$$

\vspace{0.8em}
Thus it suffices to show that the right-hand side tends to $0$ as $n,p\rightarrow +\infty$. We start by stating some new a priori estimates with our new integrability assumptions for $\xi$ and $\widehat F^0$ (see the Appendix for the proof)
\begin{align}
&\nonumber \underset{n}{\Sup}\left\{\underset{\mathbb P\in\mathcal P^\kappa_H}{\Sup}\mathbb E^\mathbb P\left[\underset{0\leq t\leq T}{\Sup}\abs{y_t^{\mathbb P,n}}^{2+\epsilon'}\right]+\underset{\mathbb P\in\mathcal P^\kappa_H}{\Sup}\mathbb E^\mathbb P\left[\left(\int_0^T\abs{\widehat a_t^{1/2}z_t^{\mathbb P,n}}^2dt\right)^{\frac{2+\epsilon'}{2}}\right]\right\}\\
&\leq C\left(1+\No{\xi}_{L^{2+\epsilon',\kappa}_H}^{2+\epsilon'}+\underset{\mathbb P\in\mathcal P^\kappa_H}{\Sup}\mathbb E^\mathbb P\left[\int_0^T\abs{\widehat F^0_s}^{2+\epsilon'}ds\right]\right).
\label{estimm}
\end{align}

\vspace{0.8em}
Let us denote $\delta\widehat{F}^{n,p}_t:=\widehat{F}^n(t,y^{\mathbb P,n}_t,z^{\mathbb P,n}_t)-\widehat{F}^p(t,y^{\mathbb P,p}_t,z^{\mathbb P,p}_t)$. Applying It\^o's formula to $|y^{\mathbb P,n}_t-y^{\mathbb P,p}_t|^2$ and taking conditional expectations yields, for all $\mathbb P \in \mathcal P_H^\kappa$
\begin{equation}
\label{abcde}
|y^{\mathbb P,n}_t-y^{\mathbb P,p}_t|^2+\mathbb E_t^\mathbb P\left[\int_t^T |\widehat{a}_s^{1/2}(z^{\mathbb P,n}_s-z^{\mathbb P,p}_s)|^2ds\right] \leq2 \mathbb{E}_t^{\mathbb{P}} \left[\int_t^T\abs{y^{\mathbb P,n}_s-y^{\mathbb P,p}_s}\abs{\delta\widehat{F}^{n,p}_s}ds\right].
\end{equation}

\vspace{0.8em}
Now since $\epsilon'>0$, it follows from Doob maximal inequality (in the classical form under a single measure) and the Cauchy-Schwarz inequality that
\begin{align}
\nonumber \mathbb E^\mathbb P\left[\underset{0\leq t\leq T}{\Sup}\scriptstyle\left(\mathbb E_t^\mathbb P\left[\int_t^T |\widehat{a}_s^{\frac12}(z^{\mathbb P,n}_s-z^{\mathbb P,p}_s)|^2ds\right]\right)^{\frac{2+\epsilon'}{2}}\right]&\leq C\mathbb E^\mathbb P\left[\left(\int_0^T\scriptstyle\abs{y^{\mathbb P,n}_s-y^{\mathbb P,p}_s}\abs{\delta\widehat{F}^{n,p}_s}ds\right)^{\frac{2+\epsilon'}{2}}\right]\\
\nonumber &\leq C\left(\mathbb E^\mathbb P\left[\int_0^T\abs{y_s^{\mathbb P,n}-y_s^{\mathbb P,p}}^{2+\epsilon'}ds\right]\right)^{\frac12}\\
&\hspace{0.9em}\times\left(\mathbb E^\mathbb P\left[\left(\int_0^T\abs{\delta\widehat{F}^{n,p}_s}^{2}ds\right)^{\frac{2+\epsilon'}{2}}\right]\right)^{\frac12}.
\label{gryui}
\end{align}

\vspace{0.9em}
By the uniform growth property $\rm{(i)}$ of Lemma \ref{lemma.approx}, we have for all $t\in[0,T]$
\begin{equation*}
|\delta\widehat{F}^{n,p}_t|^2\leq C\left(1+|\widehat F^0_t|^2+|y^{\mathbb P,n}_t|^2+|y^{\mathbb P,p}_t|^2+|\widehat{a}_t^{1/2}z^{\mathbb P,n}_t|^2+|\widehat{a}_t^{1/2}z^{\mathbb P,p}_t|^2\right), \text{ } \mathbb P-a.s., \text{ } \forall \text{ }\mathbb P \in \mathcal P_H^\kappa.
\end{equation*}

\vspace{0.8em}
Hence using the uniform a priori estimates of \reff{estimm}, we have for all $\mathbb P\in\mathcal P^\kappa_H$
\begin{align}
\label{eq.delatf}
\nonumber \mathbb{E}^{\mathbb P}\left[\left(\int_0^T|\delta\widehat{F}^{n,p}_t|^2dt\right)^{\frac{2+\epsilon'}{2}}\right]&\leq C\left(1+\underset{\mathbb P \in \mathcal P_H^\kappa}{\Sup}\mathbb{E}^{\mathbb P}\left[\left(\int_0^T\abs{\widehat F^0_s}^2ds\right)^{\frac{2+\epsilon'}{2}}\right]\right)\\
\nonumber&\hspace{0.9em}+C\left(\underset{n}{\Sup}\underset{\mathbb P \in \mathcal P_H^\kappa}{\Sup}\mathbb{E}^{\mathbb P}\left[\underset{0\leq t\leq T}{\Sup}|y_t^{\mathbb P,n}|^{2+\epsilon'}\right]\right)\\
\nonumber&\hspace{0.9em}+C\left(\underset{n}{\Sup}\underset{\mathbb P \in \mathcal P_H^\kappa}{\Sup}\mathbb{E}^{\mathbb P}\left[\left(\int_0^T|\widehat{a}_s^{1/2}z_s^{\mathbb P,n}|^2ds\right)^{\frac{2+\epsilon'}{2}}\right]\right)\\[0.8em]
&\leq C\left(1+\No{\xi}^{2+\epsilon}_{L^{2+\epsilon',\kappa}_H}+\underset{\mathbb P\in\mathcal P^\kappa_H}{\sup}\mathbb E^\mathbb P\left[\int_0^T\abs{\widehat F^0_s}^{2+\epsilon'}ds\right]\right).
\end{align}

\vspace{0.8em}
We then have
\begin{align}\label{ttt}
\underset{\mathbb P \in \mathcal P_H^\kappa}{\Sup}\mathbb E^\mathbb P\left[\int_0^T\abs{y_t^{\mathbb P,n}-y_t^{\mathbb P,p}}^{2+\epsilon'}dt\right]\leq \int_0^T\underset{\mathbb P \in \mathcal P_H^\kappa}{\Sup}\mathbb E^\mathbb P\left[\abs{y_t^{\mathbb P,n}-y_t^{\mathbb P,p}}^{2+\epsilon'}\right]dt\underset{n,p\rightarrow+\infty}{\longrightarrow}0,
\end{align}
where we used Lemma \ref{lemma.yl2} and the standard Lebesgue dominated convergence Theorem.

\vspace{0.8em}
Therefore, plugging the estimate \reff{eq.delatf} in \reff{gryui}, using \reff{ttt}, and sending $n,p$ to $+\infty$, we see that
\begin{align}
\nonumber \underset{\mathbb P\in\mathcal P^\kappa_H}{\Sup}\mathbb E^\mathbb P\left[\left(\int_0^T |\widehat{a}_s^{1/2}(z^{\mathbb P,n}_s-z^{\mathbb P,p}_s)|^2ds\right)^{\frac{2+\epsilon'}{2}}\right]&
\leq\underset{\mathbb P\in\mathcal P^\kappa_H}{\Sup}\mathbb E^\mathbb P\left[\scriptstyle\underset{0\leq t\leq T}{\Sup}\left(\mathbb E_t^\mathbb P\left[\int_t^T |\widehat{a}_s^{1/2}(z^{\mathbb P,n}_s-z^{\mathbb P,p}_s)|^2ds\right]\right)^{\frac{2+\epsilon'}{2}}\right]\\
&\underset{n,p\rightarrow+\infty}{\longrightarrow}0.
\label{convz}
\end{align}

\vspace{0.8em}
Then, by It\^o's formula we similarly obtain
\begin{align}
\nonumber \underset{\mathbb P \in \mathcal P_H^\kappa}{\Sup}\mathbb{E}^{\mathbb{P}} \left[\underset{0\leq t \leq T}{\Sup}|y^{\mathbb P,n}_t-y^{\mathbb P,p}_t|^{2+\epsilon'}\right]&\leq C\underset{\mathbb P \in \mathcal P_H^\kappa}{\Sup}\mathbb{E}^{\mathbb P}\left[\scriptstyle\underset{0\leq t \leq T}{\Sup}\abs{\int_t^T(y^{\mathbb P,n}_s-y^{\mathbb P,p}_s)(z^{\mathbb P,n}_s-z^{\mathbb P,p}_s)dB_s}^{\frac{2+\epsilon'}{2}}\right]\\
\nonumber &\hspace{0.9em}+2\left(\underset{\mathbb P \in \mathcal P_H^\kappa}{\Sup}\mathbb{E}^{\mathbb P}\left[\int_0^T|y^{\mathbb P,n}_t-y^{\mathbb P,p}_t|^{2+\epsilon'}dt\right]\right)^{1/2}\\
&\hspace{0.9em}\times\left(\underset{\mathbb P \in \mathcal P_H^\kappa}{\Sup}\mathbb{E}^{\mathbb P}\left[\left(\int_0^T|\delta\widehat{F}^{n,p}_t|^2dt\right)^{\frac{2+\epsilon'}{2}}\right]\right)^{1/2}.
\label{truc2}
\end{align}

\vspace{0.8em}
By previous calculations we know that the second term tends to $0$ as $n,p\rightarrow +\infty$. For the first one, we have by the Burkholder-Davis-Gundy inequality (where we recall that the constants involved are universal and thus do not depend on $\mathbb P$)
\begin{align*}
\mathbb{E}^{\mathbb{P}} \left[\scriptstyle  \underset{0\leq t \leq T}{\Sup} \abs{\int_t^T(y^{\mathbb{P},n}_s - y^{\mathbb{P},p}_s) (z^{\mathbb{P},n}_s - z^{\mathbb{P},p}_s) dB_s}^{\frac{2+\epsilon'}{2}}\right] & \leq C_0\mathbb{E}^{\mathbb{P}}\left[ \scriptstyle\left(\int_0^T |y^{\mathbb{P},n}_s - y^{\mathbb{P},p}_s|^2 |\widehat{a}_s^{\frac12} ( z^{\mathbb{P},n}_s - z^{\mathbb{P},p}_s )|^2ds\right)^{\frac{2+\epsilon'}{4}}\right]\\
&\leq C_0\mathbb{E}^{\mathbb P}\left[\scriptstyle\left(\underset{0\leq t\leq T}{\scriptstyle\Sup}|y^{\mathbb P,n}_t-y^{\mathbb P,p}_t|^2\int_0^T|\widehat{a}_s^{\frac12}(z^{\mathbb P,n}_s-z^{\mathbb P,p}_s)|^2ds\right)^{\frac{2+\epsilon'}{4}}\right]\\
&\leq \frac{1}{2C}\underset{\mathbb P \in \mathcal P_H^\kappa}{\Sup}\mathbb{E}^{\mathbb P}\left[\underset{0\leq t\leq T}{\Sup}|y^{\mathbb P,n}_t-y^{\mathbb P,p}_t|^{2+\epsilon'}\right]\\
&\hspace{0.9em}+\frac{CC_0^2}{4}\underset{\mathbb P \in \mathcal P_H^\kappa}{\Sup}\mathbb{E}^{\mathbb P}\left[\left(\int_0^T|\scriptstyle\widehat{a}_s^{\frac12}(z^{\mathbb P,n}_s-z^{\mathbb P,p}_s)|^2ds\right)^{\frac{2+\epsilon'}{2}}\right],
\end{align*}
where $C$ is the constant in \reff{truc2}.

\vspace{0.8em}
Reporting this in the above inequality \reff{truc2} and letting $n,p$ go to $+\infty$ then finishes the proof.
\ep

\vspace{0.8em}
\begin{Remark}
In contrast with the classical case, we proved here the convergence of $Y^n$ in $\mathbb D^{2,\kappa}_H$ before proving any convergence for $Z^n$. Proceeding in this order is crucial because of the process $K^{\mathbb P,n}$, which prevents us from using the usual techniques. Then, it is natural to use the representation formula for $Y^n$ to control the $\mathbb D^{2,\kappa}_H$ norm of $Y^n-Y^p$ by a certain norm of $y^{\mathbb P,n}-y^{\mathbb P,p}$. It turns out that we end up with a norm which is closely related to the $\mathbb L^{2,\kappa}_H$ norm. However, this norm for the process $y^{\mathbb P,n}$ is not tractable for classical BSDEs, therefore, we have to use the generalized Doob inequality (which is currently conjectured to be the best possible, see Remark $2.9$ in \cite{stz}) to return to the usual norm for $y^{\mathbb P,n}$. This in turn forces us to assume stronger integrability assumptions on $\xi$ and $\widehat F^0$.
\end{Remark}

\vspace{0.8em}
We just have proved that the sequence $(Y^n)_{n}$ is Cauchy in the Banach $\mathbb D^{2,\kappa}_H$. Thus it converges to $Y$ in $\mathbb D^{2,\kappa}_H$. Let us now focus on $Z^n$ and $K^{\mathbb P,n}$.

\begin{Lemma}\label{conv.zk}
There exist a process $Z\in\mathbb H^{2,\kappa}_H$ and a non-decreasing process $K^\mathbb P\in\mathbb D^2(\mathbb P)$ such that
\begin{equation*}
\No{Z^{n}-Z}^2_{\mathbb H_H^{2,\kappa}}+\underset{\mathbb P \in \mathcal P_H^\kappa}{\Sup}\mathbb{E}^{\mathbb P}\left[\underset{0\leq t \leq T}{\Sup}|K_t^{\mathbb P,n}-K_t^{\mathbb P}|^2\right]\longrightarrow0,\text{ as $n\rightarrow +\infty$}.
\end{equation*}
\end{Lemma}

\vspace{0.8em}
\proof
Denote $\delta\widehat{F}^{n,p}_t:=\widehat{F}_t^n(Y^{n}_t,Z^{n}_t)-\widehat{F}_t^p(Y^{p}_t,Z^{p}_t)$. Applying It\^o's formula to $|Y^{n}_t-Y^{p}_t|^2$ and taking expectations yields, for all $\mathbb P \in \mathcal P_H^\kappa$

\begin{align*}
\mathbb{E}^{\mathbb{P}} \left[ |Y^{n}_0-Y^{p}_0|^2+\int_0^T |\widehat{a}_t^{1/2}(Z^{n}_t-Z^{p}_t)|^2dt\right]& \leq2 \mathbb{E}^{\mathbb{P}} \left[\int_0^T(Y^{n}_t-Y^{p}_t)\delta\widehat{F}^{n,p}_tdt\right]\\
&\hspace{0.9em}+2 \mathbb{E}^{\mathbb{P}} \left[\int_0^T(Y^{n}_t-Y^{p}_t)d(K_t^{\mathbb P,n}-K_t^{\mathbb P,n})\right].
\end{align*}

\vspace{0.8em}
Then
\begin{align*}
\mathbb{E}^{\mathbb{P}} \left[ \int_0^T|\widehat{a}_t^{1/2}(Z^{n}_t-Z^{p}_t)|^2dt\right]&\leq
2\left(\underset{\mathbb P \in \mathcal P_H^\kappa}{\Sup}\mathbb{E}^{\mathbb P}\left[\int_0^T|Y^{n}_t-Y^{p}_t|^2dt\right]\right)^{1/2}\\
&\hspace{0.9em}\times\left(\underset{\mathbb P \in \mathcal P_H^\kappa}{\Sup}\mathbb{E}^{\mathbb P}\left[\int_0^T|\delta\widehat{F}^{n,p}_t|^2dt\right]\right)^{1/2}\\
&\hspace{0.9em} + 2\left(\underset{\mathbb P \in \mathcal P_H^\kappa}{\Sup}\mathbb{E}^{\mathbb P}\left[\underset{0\leq t \leq T}{\Sup}|Y^n_t-Y_t^p|^2\right]\right)^{1/2}\\
&\hspace{0.9em}\times\left(\underset{\mathbb P \in \mathcal P_H^\kappa}{\Sup}\mathbb{E}^{\mathbb P}\left[|K^{\mathbb P,n}_T|^2+|K^{\mathbb P,p}_T|^2\right]\right)^{1/2}.
\end{align*}

\vspace{0.8em}
Notice that the right-hand side tends to $0$ uniformly in $\mathbb P$ as $n,p\rightarrow +\infty$ due to Lemmas \ref{estim} and \ref{conv.d2y} and \reff{eq.delatf}. Thus $(Z^n)$ is a Cauchy sequence in $\mathbb H^{2,\kappa}_H$ and therefore converges to a process $Z\in \mathbb H^{2,\kappa}_H$ .

\vspace{0.8em}
Now by \reff{2bsde.k}, we have
\begin{align*}
\underset{\mathbb P \in \mathcal P_H^\kappa}{\Sup}\mathbb{E}^{\mathbb P}\left[\underset{0\leq t \leq T}{\Sup}|K_t^{\mathbb P,n}-K_t^{\mathbb P,p}|^2\right]&\leq C\underset{\mathbb P \in \mathcal P_H^\kappa}{\Sup}\mathbb{E}^{\mathbb P}\left[|Y_0^n-Y_0^p|^2+\underset{0\leq t\leq T}{\Sup}|Y_t^n-Y_t^p|^2\right]\\
&\hspace{0.9em}+C\underset{\mathbb P \in \mathcal P_H^\kappa}{\Sup}\mathbb{E}^{\mathbb P}\left[\underset{0\leq t \leq T}{\Sup}\abs{\int_t^T(Z_s^n-Z_s^p)dB_s}^2\right]\\
&\hspace{0.9em}+C\underset{\mathbb P \in \mathcal P_H^\kappa}{\Sup}\mathbb{E}^{\mathbb P}\left[\int_0^T|\delta\widehat{F}^{n,p}_t|^2dt\right].
\end{align*}

\vspace{0.9em}
The first two terms on the right-hand side tend to $0$ as $n,p\rightarrow +\infty$ thanks to Lemma \ref{conv.d2y}. For the last one, using BDG inequality and the result we just proved on the sequence $(Z^n)$, we see that it also tends to $0$. Thus, in order to finish the proof, we need to show that the term involving $\delta\widehat{F}^{n,p}$ converges to 0. This is deduced from the following facts

\begin{itemize}
\item[$\bullet$] $Y^n\nearrow Y$ in $\mathbb D^{2,\kappa}_H$ and $dt\times \mathcal P_H^\kappa-q.s$.
\item[$\bullet$] By Lemma \ref{assump.conv}, $\widehat F^n$ converges uniformly to $\widehat F$ in $(y,z)$. Moreover, we have from the Lipschitz property of $F$ in $z$ and its uniform continuity in $y$
\begin{align*}
\abs{\widehat F_t(Y_t,Z_t)-\widehat F_t^n(Y_t^n,Z_t^n)}&\leq \rho(\abs{Y_t^n-Y_t})+\abs{\widehat F_t(Y_t^n,Z_t)-\widehat F_t^n(Y_t^n,Z_t)}\\
&\hspace{0.9em}+ C|\widehat a_t^{1/2}(Z_t-Z_t^n)|,
\end{align*}
for some modulus of continuity $\rho$.

\vspace{0.8em}
When taking expectation under $\mathbb P$ and supremum over all $\mathbb P\in\mathcal P^\kappa_H$, the convergence to $0$ for the term involving $Z-Z^n$ is clear by our previous result . For the second one in the right-hand side above, we can use the same arguments as in the proof of Lemma \ref{lemma.yl2} to obtain that it also converges to $0$. Finally, we recall that since the space $\Omega$ is convex, it is a classical result that we can choose the modulus of continuity $\rho$ to be concave, non-decreasing and sub-linear. We then have by Jensen inequality
$$\underset{\mathbb P\in\mathcal P^\kappa_H}{\sup}\mathbb E^\mathbb P\left[ \rho(\abs{Y_t^n-Y_t})\right]\leq\rho\left(\underset{\mathbb P\in\mathcal P^\kappa_H}{\sup}\mathbb E^\mathbb P\left[ \abs{Y_t^n-Y_t}\right]\right)\longrightarrow 0,$$
by Lemma \ref{conv.d2y}.
\end{itemize}

\vspace{0.8em}
Consequently, for all $\mathbb P\in \mathcal P_H^\kappa$, there exists a non-decreasing and progressively measurable process $K^\mathbb P$ such that
$$\underset{\mathbb P \in \mathcal P_H^\kappa}{\Sup}\mathbb{E}^{\mathbb P}\left[\underset{0\leq t \leq T}{\Sup}|K_t^{\mathbb P,n}-K_t^{\mathbb P}|^2\right]\underset{n\rightarrow +\infty}{\rightarrow}0.$$

\vspace{0.8em}
Moreover, since the $K_t^{\mathbb P,n}$ are c\`adl\`ag, so is $K_t^{\mathbb P}.$
\ep

\vspace{1em}
{\bf Proof of Theorem \ref{th1}}. Taking limits in the 2BSDE \reff{2bsden}, we obtain that $(Y,Z)\in \mathbb D^{2,\kappa}_H\times\mathbb H^{2,\kappa}_H$ is a solution of the 2BSDE \reff{2bsde}. To conclude the proof of existence, it remains to check the minimum condition \reff{2bsde.min}. But for all $\mathbb P\in \mathcal P_H^\kappa$, we know that $K^{\mathbb P,n}$ verifies \reff{2bsde.min}. Then we can pass to the limit in the minimum condition verified by $K^{\mathbb P,n}$. Indeed, we have for all $\mathbb P$
$$K^{\mathbb P,n}_t=\underset{\mathbb P^{'}\in\mathcal P^\kappa_H(t^+,\mathbb P)}{\einf^\mathbb P}\mathbb E^{\mathbb P^{'}}\left[K^{\mathbb P^{'},n}_T\right],\text{ }\mathbb P-a.s.$$

Extracting a subsequence if necessary, the left-hand side converges to $K^{\mathbb P}_t$, $\mathcal P^\kappa_H-q.s.$ Then as in the beginning of the proof of Lemma \ref{conv.d2y}, we can write
$$\underset{\mathbb P^{'}\in\mathcal P^\kappa_H(t^+,\mathbb P)}{\einf^\mathbb P}\mathbb E^{\mathbb P^{'}}\left[K^{\mathbb P^{'},n}_T\right]=\underset{m\rightarrow +\infty}{\lim}\downarrow \mathbb E^{\mathbb P_m}\left[K^{\mathbb P_m,n}_T\right],$$
where $(\mathbb P_m)_{m\geq 0}$ is a sequence of probability measures belonging to $\mathcal P^\kappa_H(t^+,\mathbb P)$.

\vspace{0.8em}
Then by Lemma \ref{conv.zk}, we know that $\mathbb E^\mathbb P[K_T^{\mathbb P,n}]$ converges to $\mathbb E^\mathbb P[K_T^{\mathbb P}]$ uniformly in $\mathbb P\in\mathcal P^\kappa_H$. Thus we can take the limit in $n$ in the above expression and switch the limits in $n$ and $m$. This shows that
$$K_t^\mathbb P=\underset{m\rightarrow +\infty}{\lim}\downarrow \mathbb E^{\mathbb P_m}\left[K^{\mathbb P_m}_T\right]\geq \underset{\mathbb P^{'}\in\mathcal P^\kappa_H(t^+,\mathbb P)}{\einf^\mathbb P}\mathbb E^{\mathbb P^{'}}\left[K^{\mathbb P^{'}}_T\right].$$

The converse inequality is trivial since the process $K^\mathbb P$ is non-decreasing. Thus, the minimum condition is satisfied, and the proof of Theorem \ref{th1} is complete.
\ep

\begin{Remark}
In comparison with the classical BSDE framework, we had to add some assumptions here to prove existence of a solution. The question is whether these assumptions can be weakened by using another construction for the solution. For instance, we may use the so called regular conditional probability distributions as in \cite{stz} and \cite{posz}. However, as mentioned in Remark $4.1$ in \cite{posz}, even though we could construct a candidate solution when the terminal condition is in $UC_b(\Omega)$, when trying to check that the family $\left\{K^\mathbb P,\mathbb P\in\mathcal P^\kappa_H\right\}$ obtained verifies the minimum condition, our monotonicity assumption is not sufficient. Thus, regardless of the solution construction method, we have to add some assumptions to prove existence.
\end{Remark}

\vspace{0.8em}
\paragraph{Acknowledgements} The author is very grateful to Nizar Touzi and Chao Zhou for their precious suggestions throughout his work. The author would also like to thank two anonymous referees and an associate editor for their helpful remarks and comments which improved an earlier version of this paper.

\vspace{0.8em}
\begin{appendix}
\section{Appendix}
\begin{Proposition}\label{doob}
Let $\xi$ be some $\mathcal F_T$ random variable, and let $n\geq 1$. Then, for all $p>n$ 
$$\underset{\mathbb P\in\mathcal P^\kappa_H}{\Sup}\mathbb E^\mathbb P\left[\underset{0\leq t \leq T}{\Sup}\left(\underset{\mathbb P^{'}\in\mathcal P^\kappa_H(t^+,\mathbb P)}{\esup^\mathbb P}\mathbb E_t^{\mathbb P^{'}}\left[|\xi|\right]\right)^n\right]\leq C\underset{\mathbb P\in\mathcal P^\kappa_H}{\Sup}\left(\mathbb E^\mathbb P\left[\abs{\xi}^{p}\right]\right)^{\frac{n}{p}}.$$

\end{Proposition}

\vspace{0.8em}
\proof
The proof follows the ideas of the proof of Lemma $6.2$ in \cite{stz3}, which closely follows the classical proof of the Doob maximal inequality (see also \cite{song} for related results).

\vspace{0.8em}
Fix some $\mathbb P$, and let us note $$X_t^\mathbb P:=\underset{\mathbb P^{'}\in\mathcal P^\kappa_H(t^+,\mathbb P)}{\esup^\mathbb P}\mathbb E_t^{\mathbb P^{'}}\left[|\xi|\right].$$

\vspace{0.8em}
Then, it can be shown that $X_t^\mathbb P$ is a $\mathbb P$-supermartingale, and thus admits a c\`adl\`ag version. For all $\lambda>0$, let us define the following $\mathbb F^+$ stopping times
$$\tau^\mathbb P_\lambda=\inf\left\{t\geq 0, \text{ }X_t^\mathbb P\geq \lambda, \text{ }\mathbb P-a.s.\right\}.$$

\vspace{0.8em}
Define $X^{\mathbb P,*}:=\underset{0\leq t \leq T}{\Sup}X_t^\mathbb P$. Then, we have
$$\mathbb P(X^{\mathbb P,*}\geq \lambda)=\mathbb P(\tau^\mathbb P_\lambda\leq T)\leq \frac1\lambda\mathbb E^\mathbb P\left[X^\mathbb P_{\tau^\mathbb P_\lambda}1_{\tau^\mathbb P_\lambda\leq T}\right].$$

\vspace{0.8em}
As previously in this chapter, we know that there exists a sequence $(\mathbb P_n)_{n\geq 0}\subset\mathcal P^\kappa_H\left((\tau^\mathbb P_\lambda)^+,\mathbb P\right)$ such that
$$X^\mathbb P_{\tau^\mathbb P_\lambda}=\underset{n\rightarrow+\infty}{\lim}\uparrow\mathbb E_{\tau^\mathbb P_\lambda}^{\mathbb P_n}\left[\abs{\xi}\right].$$

\vspace{0.8em}
Hence, using in this order the monotone convergence theorem and the fact that all the $\mathbb P_n$ coincide with $\mathbb P$ on $\mathcal F_{\tau_\lambda^\mathbb P}$, we have

\begin{align*}
\mathbb E^\mathbb P\left[X^\mathbb P_{\tau^\mathbb P_\lambda}1_{\tau^\mathbb P_\lambda\leq T}\right]&=\underset{n\rightarrow+\infty}{\lim}\uparrow\mathbb E^\mathbb P\left[\mathbb E_{\tau^\mathbb P_\lambda}^{\mathbb P_n}\left[\abs{\xi}\right]1_{\tau^\mathbb P_\lambda\leq T}\right]\\
&=\underset{n\rightarrow+\infty}{\lim}\uparrow\mathbb E^{\mathbb P_n}\left[\mathbb E_{\tau^\mathbb P_\lambda}^{\mathbb P_n}\left[\abs{\xi}\right]1_{\tau^\mathbb P_\lambda\leq T}\right]\\
&=\underset{n\rightarrow+\infty}{\lim}\uparrow\mathbb E^{\mathbb P_n}\left[\abs{\xi}1_{\tau^\mathbb P_\lambda\leq T}\right]\\
&\leq\underset{n\rightarrow+\infty}{\lim}\uparrow\left(\mathbb E^{\mathbb P_n}\left[\abs{\xi}^p\right]\right)^{1/p}\mathbb P_n(\tau^\mathbb P_\lambda\leq T)^{1-\frac1p}\\
&\leq\mathbb P(\tau^\mathbb P_\lambda\leq T)^{1-\frac1p}\underset{n\rightarrow+\infty}{\lim}\uparrow\left(\mathbb E^{\mathbb P_n}\left[\abs{\xi}^p\right]\right)^{1/p}\\
&\leq \mathbb P(\tau^\mathbb P_\lambda\leq T)^{1-\frac1p}\underset{\mathbb P\in\mathcal P^\kappa_H}{\Sup}\left(\mathbb E^\mathbb P\left[\abs{\xi}^p\right]\right)^{1/p}.
\end{align*}

\vspace{0.8em}
Using this estimate, we finally get
$$\mathbb P(X^{\mathbb P,*}\geq \lambda)\leq \frac{1}{\lambda^p}\underset{\mathbb P\in\mathcal P^\kappa_H}{\Sup}\mathbb E^\mathbb P\left[\abs{\xi}^p\right].$$

\vspace{0.8em}
Now, for every $\lambda_0>0$, we have
\begin{align*}
\mathbb E^\mathbb P\left[(X^{\mathbb P,*})^n\right]&=n\int_0^{+\infty}\lambda^{n-1}\mathbb P(X^{\mathbb P,*}\geq \lambda)d\lambda\\
&=n\int_0^{\lambda_0}\lambda^{n-1}\mathbb P(X^{\mathbb P,*}\geq \lambda)d\lambda+n\int_{\lambda_0}^{+\infty}\lambda^{n-1}\mathbb P(X^{\mathbb P,*}\geq \lambda)d\lambda\\
&\leq \lambda_0^n+n\underset{\mathbb P\in\mathcal P^\kappa_H}{\Sup}\mathbb E^\mathbb P\left[\abs{\xi}^p\right]\int_{\lambda_0}^{+\infty}\frac{d\lambda}{\lambda^{p-n+1}}\\
&=\lambda_0^n+n\underset{\mathbb P\in\mathcal P^\kappa_H}{\Sup}\mathbb E^\mathbb P\left[\abs{\xi}^p\right]\frac{\lambda_0^{n-p}}{p-n}.
\end{align*}

\vspace{0.8em}
Choosing $\lambda_0=\underset{\mathbb P\in\mathcal P^\kappa_H}{\Sup}\left(\mathbb E^\mathbb P\left[\abs{\xi}^p\right]\right)^{1/p}$, this yields
$$\mathbb E^\mathbb P\left[(X^{\mathbb P,*})^n\right]\leq \left(1+\frac{n}{p-n}\right)\underset{\mathbb P\in\mathcal P^\kappa_H}{\Sup}\left(\mathbb E^\mathbb P\left[\abs{\xi}^p\right]\right)^{n/p}.$$

\vspace{0.8em}
Hence the result.
\ep

\vspace{0.8em}
{\bf Proof of Lemma \ref{assump.conv}}. The uniform convergence result is a simple consequence of a result already proved by Lasry and Lyons in \cite{ll} (see the Theorem on page $4$ and Remark $\rm{(iv)}$ on page $5$). For the inequality, since $\widehat F$ is uniformly continuous in $y$, there exists a modulus of continuity $\rho$ with linear growth. Then, it follows that

\begin{align*}
F_t-F^n_t&=\underset{u\in \mathbb Q,u\geq y}{\sup}\left\{F_t(y,z,a)-F_t(u,z,a)-n\abs{y-u}\right\}\\
&\leq \underset{u\in \mathbb Q,u\geq y}{\sup}\left\{C\rho(\abs{y-u})-n\abs{y-u}\right\}\\
&=\underset{u\geq 0}{\sup}\left\{C\rho(u)-nu\right\}.
\end{align*} 

\vspace{0.8em}
Since $\rho$ has linear growth in $y$, the function on the right-hand side above is clearly decreasing in $n$ and thus is dominated by a constant, which gives us the result.

\vspace{0.8em}
Finally, the equality of the two suprema is clear because the quantity considered here is bounded and continuous in $\omega$.
\ep

\vspace{0.8em}
\begin{Proposition}
Assume that there exists an $\epsilon>0$ such that $\xi\in L^{2+\epsilon,\kappa}_H$ and 
$$\underset{\mathbb P\in\mathcal P^\kappa_H}{\Sup}\mathbb E^\mathbb P\left[\int_0^T\abs{\widehat F^0_t}^{2+\epsilon}dt\right]<+\infty.$$

\vspace{0.8em}
Then we have
\begin{align*}
&\nonumber \underset{n}{\Sup}\left\{\underset{\mathbb P\in\mathcal P^\kappa_H}{\Sup}\mathbb E^\mathbb P\left[\underset{0\leq t\leq T}{\Sup}\abs{y_t^{\mathbb P,n}}^{2+\epsilon}\right]+\underset{\mathbb P\in\mathcal P^\kappa_H}{\Sup}\mathbb E^\mathbb P\left[\left(\int_0^T\abs{\widehat a_t^{1/2}z_t^{\mathbb P,n}}^2dt\right)^{\frac{2+\epsilon}{2}}\right]\right\}\\
&\leq C\left(1+\No{\xi}_{L^{2+\epsilon,\kappa}_H}^{2+\epsilon}+\underset{\mathbb P\in\mathcal P^\kappa_H}{\Sup}\mathbb E^\mathbb P\left[\int_0^T\abs{\widehat F^0_s}^{2+\epsilon}ds\right]\right).
\end{align*}
\end{Proposition}

\vspace{0.8em}
\proof
Fix a $\mathbb P$, and consider the BSDE \reff{bsde}. As in the proof of Lemma \ref{estim}, we have that for all $n$, $y^{\mathbb P,n}\leq u^\mathbb P$, $\mathbb P$-a.s. Now let $\alpha$ be some positive constant which will be fixed later and let $\eta\in(0,1)$. By It\^o's formula we have

\begin{align*}
e^{\alpha t}\abs{u_t^\mathbb P}^2+\int_t^Te^{\alpha s}\abs{\widehat a_s^{1/2}v_s^\mathbb P}^2ds&=e^{\alpha T}\abs{\xi}^2+2\int_t^Te^{\alpha s}u_s^\mathbb P\left(\abs{\widehat F^0_s}+C\right)ds\\
&\hspace{0.9em}+2C\int_t^Tu_s^\mathbb P\left(e^{\alpha s}\abs{u_s^\mathbb P}+\abs{\widehat a_s^{1/2}v_s^\mathbb P}\right)ds-\alpha\int_t^Te^{\alpha s}\abs{u_s^\mathbb P}^2ds\\
&\hspace{0.9em}-2\int_t^Te^{\alpha s}u_s^\mathbb Pv_s^\mathbb PdB_s\\
&\leq e^{\alpha T}\abs{\xi}^2+2\int_t^Te^{\alpha s}\abs{u_s^\mathbb P}\left(\abs{\widehat F^0_s}+C\right)ds\\
&\hspace{0.9em}+\left(2C+\frac{C^2}{\eta}-\alpha\right)\int_t^Te^{\alpha s}\abs{u_s^\mathbb P}^2ds +\eta\int_t^Te^{\alpha s}\abs{\widehat a_s^{\frac12}v_s^\mathbb P}^2ds\\
&\hspace{0.9em}-2\int_t^Te^{\alpha s}u_s^\mathbb Pv_s^\mathbb PdB_s.
\end{align*}

\vspace{0.8em}
Now choose $\alpha$ such that $\nu:=\alpha -2C-\frac{C^2}{\eta}\geq 0$. We obtain
\begin{align}
\nonumber e^{\alpha t}\abs{u_t^\mathbb P}^2+(1-\eta)\int_t^Te^{\alpha s}\abs{\widehat a_s^{1/2}v_s^\mathbb P}^2ds+\nu\int_t^Te^{\alpha s}\abs{u_t^\mathbb P}^2ds&\leq e^{\alpha T}\abs{\xi}^2+2\int_t^Te^{\alpha s}\abs{u_s^\mathbb P\widehat F^0_s}ds\\
\nonumber&\hspace{0.9em}+2C\int_t^Te^{\alpha s}\abs{u_s^\mathbb P}ds\\
&\hspace{0.9em}-2\int_t^Te^{\alpha s}u_s^\mathbb Pv_s^\mathbb PdB_s.
\label{bidule}
\end{align}

\vspace{0.8em}
Taking conditional expectation in \reff{bidule} yields
$$e^{\alpha t}\abs{u_s^\mathbb P}^2\leq\mathbb E^\mathbb P_t\left[e^{\alpha T}\abs{\xi}^2+2\int_t^Te^{\alpha s}\abs{u_s^\mathbb P}\left(\abs{\widehat F^0_s}+C\right)ds\right].$$

\vspace{0.8em}
By Doob's maximal inequality, we then get for all $\beta\in(0,1)$, since $\epsilon>0$
\begin{align*}
\mathbb E^\mathbb P\left[\underset{0\leq t\leq T}{\Sup}e^{\alpha t\frac{2+\epsilon}{2}}\abs{u_t^\mathbb P}^{2+\epsilon}\right]&\leq C\mathbb E^\mathbb P\left[e^{\alpha T\frac{2+\epsilon}{2}}\abs{\xi}^{2+\epsilon}\right]\\
&\hspace{0.9em}+C\mathbb E^\mathbb P\left[\underset{0\leq t\leq T}{\Sup}\left(e^{\alpha t\frac{2+\epsilon}{4}}\abs{u_s^\mathbb P}^{\frac{2+\epsilon}{2}}\right)\left(\int_0^TC+\abs{\widehat F^0_s}ds\right)^{\frac{2+\epsilon}{2}}\right]\\
&\leq C\left(1+\No{\xi}_{L^{2+\epsilon,\kappa}_H}^{2+\epsilon}\right)+\beta\mathbb E^\mathbb P\left[\underset{0\leq t\leq T}{\Sup}e^{\alpha t\frac{2+\epsilon}{2}}\abs{u_t^\mathbb P}^{2+\epsilon}\right]\\
&\hspace{0.9em}+\beta \frac{C^2}{4}\mathbb E^\mathbb P\left[\int_0^T\abs{\widehat F^0_s}^{2+\epsilon}ds\right].
\end{align*}

\vspace{0.8em}
Since $\alpha$ is positive and $y^{\mathbb P,n}\leq u^\mathbb P$, we get finally for all $n$
\begin{equation}
\label{estim.u}
\underset{\mathbb P\in\mathcal P^\kappa_H}{\Sup}\mathbb E^\mathbb P\left[\underset{0\leq t\leq T}{\Sup}\abs{y_t^{\mathbb P,n}}^{2+\epsilon}\right]\leq C\left(1+\No{\xi}_{L^{2+\epsilon,\kappa}_H}^{2+\epsilon}+\underset{\mathbb P\in\mathcal P^\kappa_H}{\Sup}\mathbb E^\mathbb P\left[\int_0^T\abs{\widehat F^0_s}^{2+\epsilon}ds\right]\right).
\end{equation}

\vspace{0.8em}
Now apply It\^o's formula to $\abs{y^{\mathbb P,n}}^2$ and put each side to the power $\frac{2+\epsilon}{2}$, we have easily
\begin{align*}
\left(\int_0^T\abs{\widehat a_t^{1/2}z_t^{\mathbb P,n}}^2dt\right)^{\frac{2+\epsilon}{2}}&\leq C\left(\abs{\xi}^{2+\epsilon}+\left(\int_0^T\abs{y_t^{\mathbb P,n}\widehat F^n_t(y_t^{\mathbb P,n},z_t^{\mathbb P,n})}dt\right)^{\frac{2+\epsilon}{2}}\right)\\
&\hspace{0.9em}+C\left(\int_0^Ty^{\mathbb P,n}_tz^{\mathbb P,n}_tdB_t\right)^{\frac{2+\epsilon}{2}}
\end{align*} 

\vspace{0.8em}
Since $\widehat F^n$ satisfies a uniform linear growth property, we get
\begin{align*}
\mathbb E^\mathbb P\left[\left(\int_0^T\abs{\widehat a_t^{1/2}z_t^{\mathbb P,n}}^2dt\right)^{\frac{2+\epsilon}{2}}\right]&\leq C\No{\xi}_{L^{2+\epsilon,\kappa}_H}^{2+\epsilon}+C\left(1+\underset{\mathbb P\in\mathcal P^\kappa_H}{\Sup}\mathbb E^\mathbb P\left[\int_0^T\abs{\widehat F^0_t}^{2+\epsilon}dt \right]\right)\\
&\hspace{0.9em}+C\underset{\mathbb P\in\mathcal P^\kappa_H}{\Sup}\mathbb E^\mathbb P\left[\underset{0\leq t\leq T}{\sup}\abs{y_t^{\mathbb P,n}}^{2+\epsilon}\right]\\
&\hspace{0.9em}+\frac13\mathbb E^\mathbb P\left[\left(\int_0^T\abs{\widehat a_t^{1/2}z_t^{\mathbb P,n}}^2dt\right)^{\frac{2+\epsilon}{2}}\right]\\
&\hspace{0.9em}+C\mathbb E^\mathbb P\left[\left(\int_0^Ty^{\mathbb P,n}_tz^{\mathbb P,n}_tdB_t\right)^{\frac{2+\epsilon}{2}}\right]
\end{align*}

\vspace{0.8em}
Hence, we get for all $\gamma\in (0,1)$ by BDG inequality and \reff{estim.u}
\begin{align*}
\mathbb E^\mathbb P\left[\left(\int_0^T\abs{\widehat a_t^{1/2}z_t^{\mathbb P,n}}^2dt\right)^{\frac{2+\epsilon}{2}}\right]&\leq C\left(1+\No{\xi}_{L^{2+\epsilon,\kappa}_H}^{2+\epsilon}+\underset{\mathbb P\in\mathcal P^\kappa_H}{\Sup}\mathbb E^\mathbb P\left[\int_0^T\abs{\widehat F^0_t}^{2+\epsilon}dt \right]\right)\\
&\hspace{0.9em}+C\mathbb E^\mathbb P\left[\left(\int_0^T\abs{y^{\mathbb P,n}_t}^2\abs{\widehat a_t^{1/2}z^{\mathbb P,n}_t}^2dt\right)^{\frac{2+\epsilon}{4}}\right]\\
& \leq C\left(1+\No{\xi}_{L^{2+\epsilon,\kappa}_H}^{2+\epsilon}+\underset{\mathbb P\in\mathcal P^\kappa_H}{\Sup}\mathbb E^\mathbb P\left[\int_0^T\abs{\widehat F^0_t}^{2+\epsilon}dt \right]\right)\\
&\hspace{0.9em}+\gamma \mathbb E^\mathbb P\left[\left(\int_0^T\abs{\widehat a_t^{1/2}z^{\mathbb P,n}_t}^2dt\right)^{\frac{2+\epsilon}{2}}\right]\\
&\hspace{0.9em}+ \frac{C^2}{4\gamma}\mathbb E^\mathbb P\left[\underset{0\leq t\leq T}{\Sup}\abs{y_t^{\mathbb P,n}}^{2+\epsilon}\right],
\end{align*}
which ends the proof.
\ep

\end{appendix}
\end{document}